\input amstex

\documentstyle{amsppt}
\loadbold

 \NoRunningHeads
 \NoBlackBoxes

\topmatter
\title
Krieger's type \\for ergodic nonsingular Poisson actions of non-(T) locally compact groups
\endtitle

\author  Alexandre I. Danilenko 
\endauthor

\abstract It is shown that each locally compact second countable  non-(T) group $G$ admits  non-strongly ergodic weakly mixing IDPFT Poisson actions of any possible Krieger type. 
These actions are amenable if and only if $G$ is amenable.
If $G$ has the Haagerup property then (and only then) 
these actions
can be chosen  of 0-type.
If $G$ is amenable  then $G$ admits weakly mixing Bernoulli actions of arbitrary Krieger type.
\endabstract

\address
B. Verkin Institute for Low Temperature Physics \& Engineering
of Ukrainian National Academy of Sciences,
47 Nauky Ave.,
 Kharkiv, 61164, UKRAINE
\endaddress
\email            alexandre.danilenko\@gmail.com
\endemail

\endtopmatter

\document

\head 0. Introduction
\endhead

This article continues investigation of ergodic properties of nonsingular Poisson suspensions initiated by Z.~Kosloff, E.~Roy and the present author in a series of works \cite{DaKoRo1}, \cite{DaKoRo2} and \cite{DaKo}.
In this paper  we deal with the class of  locally compact second countable groups without Kazhdan property~(T) and  a subclass of these groups with  the Haagerup property (or a-$T$-menable groups according to M.~Gromov).
For properties of  these groups and their applications  in group representations theory, geometric group theory, operator algebras, ergodic theory,  the Baum-Connes conjecture etc.,  we refer to the books \cite{BeHaVa} and \cite{Ch--Va} and references therein.


Our purpose is to find all possible Krieger types for the ergodic conservative Poisson actions 
of groups from the aforementioned  two classes.
We recall some standard definitions.
Let $G$ be a locally compact second countable group.
\roster
\item"---" A weakly continuous unitary representation $V=(V(g))_{g\in G}$ of $G$ in a separable Hilbert space $\Cal H$
 has {\it almost invariant vectors} if for each $\epsilon>0$, and every compact subset $K\subset G$, there is a unit vector $\xi\in\Cal H$ such that $\sup_{g\in K}\|V(g)\xi-\xi\|<\epsilon$.
\item"---" $G$ has {\it Kazhdan's property (T)} if for each unitary representation of $G$ that has almost invariant vectors,  there is  a non-zero invariant vector.
\item"---" $G$ has {\it the Haagerup property} if there is a unitary representation  $V$ of $G$ with almost invariant vectors and such that 
$V(g)\to 0$ as $g\to\infty$ in the weak operator topology. 
\endroster
We note that each amenable group has the Haagerup property (consider the left regular representation of $G$) and no group with the Haagerup property has  property~(T).

As it turned out, our main results provide new dynamical criteria for the property~(T) and the Haagerup property (see Theorems~A and B below). 
Prior to state them, we recall the famous dynamical criterion for the Kazhdan property due to A.~Connes and B.~Weiss: $G$ is non-(T) if and only if there is a weakly mixing probability preserving $G$-action that is not strongly ergodic \cite{CoWe}.
Another criterion, in terms of infinite measure preserving actions, is provided in \cite{Jo}:
$G$
is non-(T) if and only if there is  an infinite measure preserving $G$-action whose 
Koopman representation is weakly mixing and admits almost invariant vectors.
Similar criteria were proved for 
 the Haagerup property.
  The following are equivalent:
\roster
\item"---" $G$ has the Haagerup property.
\item"---" there is a mixing probability preserving $G$-action that is not strongly ergodic \cite{Ch--Va}.
\item"---" there is a 0-type infinite measure preserving $G$-action 
 whose 
Koopman representation  admits almost invariant vectors \cite{DeJoZu}.
\endroster
We now characterize the (T)-property and the Haagerup property in terms of nonsingular Poisson actions.

\proclaim{Theorem A} Let 
$\Cal K:= \{III_\lambda\mid 0\le\lambda\le1\}\cup\{II_\infty\}$.
 The following are equivalent.
\roster
\item  $G$ does not have  Kazhdan's property (T).
\item 
Given $K\in\Cal K$, 
 there is an ergodic nonsingular Poisson  $G$-action  of Krieger type $K$.
 \item
  For each $K\in\Cal K$, there is 
  an ergodic nonsingular Poisson  $G$-action $T^*$ of Krieger type $K$
 with the following extra properties:
 $T^*$  is (essentially) free,  of infinite ergodic index, non-strongly ergodic and IDPFT.
\endroster
\endproclaim

\proclaim{Theorem B} 
Let 
$\Cal K$ be as in Theorem~A.
 The following are equivalent.
\roster
\item  $G$ is a non-compact group with the Haagerup property.
\item 
There is an ergodic nonsingular Poisson  $G$-action  that is of 0-type 
 and non-strongly ergodic.
\item
For each 
 $K\in \Cal K$,
 there is an ergodic nonsingular  Poisson  $G$-action $T^*$ of Krieger type $K$ with the following extra properties:
 $T^*$  is  free,  of infinite ergodic index, of 0-type, non-strongly ergodic and IDPFT.
\endroster
\endproclaim

All necessary definitions can be found  in Section~1.
In particular, the (less common) concept of IDPFT is given in Definition~1.15.

\remark{Remark C} Every IDPFT action is amenable in the sense of Greenleaf.
Hence each Poisson $G$-action $T^*$ that appears  in Theorem~A(3) or Theorem~B(3)
is amenable  in the sense of Greenleaf.
This, in turn, implies  that  $T^*$   is amenable  (in the sense of Zimmer) if and only if $G$ is amenable.

 \endremark

We note that Theorem~A refines \cite{DaKoRo1, Theorem~8.1}, where it was shown that  $G$ has property (T) if and only if each nonsingular Poisson $G$-action admits an absolutely continuous invariant probability measure.

 \comment
 
 \remark{Remark D} While proving Theorems~A and B we also
 show 
that if $T^*$ is of type $II_\infty$ than 
there exists a $T^*$-F{\o}lner sequence of subsets of measure 1 with respect to an equivalent invariant infinite $\sigma$-finite measure.
\endremark

\endcomment


\comment

Theorems A and B significantly extend and refine the following  recent results.
\roster
\item"---" There exist  Poisson $G$-actions that  do not admit an equivalent invariant probability measure \cite{DaKoRo1, Theorem~8.2}.
\item"---" There exist nonsingular $G$-actions which are amenable in the Greenleaf sense and 
do not admit an equivalent invariant probability measure \cite{ArIsMa, Corollary~6.5}.
\item"---" There exist weakly mixing IDPFT Poisson actions of all Krieger types
 for each infinite countable discrete amenable group  \cite{DaKo}.
\item"---" 
There exist weakly mixing infinite measure preserving IDPFT $G$-actions admitting nontrivial F{\o}lner sequences \cite{Da, Theorem~D(ii)} (and partly in \cite{Jo}).
\item"---" If $G$ has the Haagerup property then there exists a weakly mixing 0-type infinite measure preserving IDPFT $G$-action admitting nontrivial F{\o}lner sequences \cite{Da, Theorem~A} (and partly in \cite{DeJoZu}).
\item"---" If $G$ has the Haagerup property then there exists a weakly mixing 0-type nonsingular $G$-action which is amenable in the Greenleaf sense and of Krie\-ger's type $III_1$ \cite{ArIsMa, Theorem~7.7}.
\endroster

\endcomment


We now state a couple of applications of Theorems~A and B.
First, as far as we know, the following problem is still open:
\roster
\item"---" {\sl Given a non-amenable $G$, is there an ergodic (or weakly mixing) non-amenable free nonsingular $G$-action $T$ that does not admit an equivalent invariant probability measure? What is the  Krieger type of $T$}?
\endroster
A partial solution was obtained recently  in \cite{ArIsMa},  where ergodic non-amenable nonsingular  $G$-actions of type $III_1$ were constructed for each non-(T) group $G$.
We generalize this result: 
it follows from  Theorem~A and Remark~C that  for each non-(T) group $G$,
 there are weakly mixing non-amenable nonsingular Poisson $G$-actions of each
 possible Krieger type $K\in\Cal K$. 


Secondly, we introduce the concept of nonsingular Bernoulli action for an arbitrary non-compact locally compact second countable group $G$.
Only the case  of discrete $G$ was considered in the literature so far.
By analogy with the probability preserving case (studied in  \cite{OrWe}), we define nonsingular Bernoulli actions as nonsingular Poisson suspensions of nonsingular totally dissipative $G$-actions.
Then we
 obtain as a corollary from Theorem~B the following result.

\proclaim{Theorem D} Let $G$ be an amenable non-compact locally compact second countable group.
If $K\in\{III_\lambda\mid 0\le \lambda\le 1\}\cup\{II_\infty\}$  
then there is a  free nonsingular IDPFT Bernoulli $G$-action of infinite ergodic index
and of  Krieger type $K$.
\endproclaim

In the case of discrete $G$, Theorem~D
 was proved in \cite{DaKo} and  \cite{BeVa} (see also \cite{VaWa} for $K=III_1$ and  \cite{KoSo} for $K=III_\lambda$ with $0<\lambda<1$).

Thirdly, we discuss an interplay between nonsingular Poisson and nonsingular Gaussian actions.
 While each nonsingular $G$-action  generates a unitary Koopman representation of $G$, each nonsingular Poisson $G$-action $T$ generates also a certain   affine Koopman representation of $G$ \cite{DaKoRo1}.
The latter gives rise to a nonsingular Gaussian $G$-action $S$ that is unitarily equivalent to $T$ \cite{DaLe, Remark~3.6}.
We recall that the nonsingular Gaussian $G$-actions were
 introduced in \cite{ArIsMa} (see also \cite{DaLe, \S3} for an alternative exposition).
A natural problem  arises: to compare  non-spectral dynamical properties of $T$ and $S$.
We show that these properties can be quite different
(see also~Proposition~8.1 for a slightly more general result).

 \proclaim{Proposition E} If $G$ is non-amenable group with the Haagerup property 
 and $T^*$ is a Poisson $G$-action (of any Krieger type) from Theorem~B(3) then the corresponding Gaussian $G$-action
 $S$ is weakly mixing,  0-type, amenable in the Greenleaf sense 
  and of Krieger type $III_1$.
 \endproclaim
 
 Thus, while there are  similarities between the theory of  nonsingular Poisson actions and the theory of nonsingular  Gaussian actions\footnote{Compare \cite{DaKoRo1} and \cite{DaKoRo2} on one side with \cite{ArIsMa}, \cite{DaLe} and \cite{MaVa} on the other side.}, Proposition~E illustrates that the nonsingular Poisson actions have richer orbital structure and they are presumably more suitable for applications in von Neumann algebras\footnote{Currently, the all known ergodic conservative nonsingular Gaussian $G$-actions are either of type $II_1$ or of type $III_1$ (see \cite{ArIsMa}, \cite{DaLe} and \cite{MaVa}).}.

As for the proof  of Theorems~A and B, we note that
 in both cases, the implications $(3)\Rightarrow(2)\Rightarrow (1)$ are straightforward.
 The main results of the paper are to prove $(1)\Rightarrow(3)$.
 Thus, let $G$ be a non-(T) group.
 We have to construct  Poisson $G$-actions satisfying the properties listed in (3).
The idea of the construction is as follows.
 Let $S=(S_g)_{g\in G}$ be a measure preserving $G$-action on an infinite $\sigma$-finite measure space $(Y,\kappa)$.
Fix a sequence $(F_n)_{n=1}^\infty$ of subsets of $Y$ such that
$\kappa(F_n)=1$ for each $n\in\Bbb N$ and $\kappa(S_gF_n\triangle F_n)\to 0$
uniformly on the compacts in $G$.
The existence of $(Y,\kappa,S,(F_n)_{n=1}^\infty)$ was proved in \cite{Jo} (see also \cite{Da}).
 Let $(X,\nu, T)$ denote the infinite sum $\bigsqcup_{n=1}^\infty(Y,\nu_n,S)$ of countably many copies of $(Y,\kappa, S)$, where $\nu_n=\kappa$ for each $n\in\Bbb N$.
The Poisson suspension $(X^*,\nu^*, T^*)$ of $(X,\nu, T)$ is canonically isomorphic to the infinite direct product  $\bigotimes_{n=1}^\infty(Y^*,\nu_n^*, S^*)$  of Poisson suspensions of $(Y,\nu_n, S)$.
 We note that $T^*$ preserves $\nu^*$.
 Now, we replace $\nu_n$ with an equivalent measure $\mu_n$ in such a way that 
 \roster
 \item"---" $\nu_n\restriction(Y\setminus F_n)=\mu_n\restriction(Y\setminus F_n)$ and $\nu_n^*\sim\mu_n^*$ for each $n$,
\item"---" the product measure $\mu^*:=\bigotimes_{n\in\Bbb N}\mu_n^*$ is  $T^*$-nonsingular and 
\item"---" $\mu^*$ is mutually singular with respect to $\bigotimes_{n\in\Bbb N}\nu_n^*$.
\endroster
 All the nonsingular Poisson actions 
 from 
 Theorems~A(3) and  B(3) appear as $(X^*,\mu^*, T^*)$ for appropriately  chosen ``parameters''  $(Y,\kappa,S, (F_n)_{n=1}^\infty, (\mu_n)_{n=1}^\infty)$.
Hence, all these Poisson actions have the 
 IDPFT structure:
 $\mu^*$ splits into an infinite direct product of  measures, each of which  admits an equivalent   invariant probability. 
We have to show that $T^*$ is conservative and  ergodic and compute the Krieger type of $T^*$.
The other properties are  more or less straightforward.
First, the conservativeness of $T^*$ is established via a (modified)  criterion from \cite{DaKoRo1}.
Then we note that under certain conditions,  each
conservative IDPFT  action is ergodic \cite{DaLe}.
In the general case considered in Theorem~A,  these conditions are not satisfied.
However they are satisfied partly,   ``along a direction''.
Therefore we  prove ``directional'' counterparts of the criterion from \cite{DaLe} as auxiliary results.
As far as we know, this ``directional'' approach originates from \cite{ArIsMa}, where it was used to show weak
mixing of certain nonsingular Gaussian actions. 
Finally, to compute the Krieger type of $T^*$,  we use again the ``directional'' approach,  the IDPFT property plus some additional tools: 
\roster
\item"---" restricted infinite products of probability preserving systems \cite{Hi} and approximation techniques from \cite{Da} (in the case of  type $II_\infty$), 
\item"---" Moore-Hill products \cite{Hi} (in the case  of  type  $III_0$),
\item"---" asymptotic properties of Skellam distributions (in the case 
of  type $III_\lambda$, $0<\lambda\le 1$).
\endroster

 The outline of the paper is as follows.
  Section~1 contains preliminary results  on unitary representations,  mixing properties and strong ergodicity of nonsingular actions, Krieger type, Maharam extensions and associated flows, Moore-Hill restricted infinite products of probability measures, IDPFT-actions,  nonsingular Poisson suspensions and amenable actions.
 Some of these results  are well known but some are new. 
 We provide complete proof in the latter case.

 In Section~2 we present  a general construction the Poisson suspension $(X^*,\mu^*,T^*)$
 depending on the parameters $(Y,\kappa,S, (F_n)_{n=1}^\infty, (\mu_n)_{n=1}^\infty)$.
Some conditions on the parameters  are found under which $(X^*,\mu^*,T^*)$ is well defined, non-strongly ergodic, etc. 

In Sections 3, 4, 5  and 6 we prove Theorems~A and B in the case $K=II_\infty$,
 $K=III_0$,  $K=III_\lambda$ with $0<\lambda<1$ and $K=III_1$ respectively.
 
 \comment
 The scheme of the proof is similar to that of the main result from \cite{DaKo}, where the case of countable amenable $G$ was considered.
 The desired Poisson actions appeared there as the Poisson suspensions of countable direct sums of the $G$-action on itself by left rotations; and the corresponding quasi-invariant measures were the counting (i.e. Haar) measures perturbed in a certain way along a fixed F{\o}lner sequence in $G$.
 In the framework of the present paper, $G$ is an arbitrary locally compact second countable group without property (T).
 Thus, $G$ is not necessarily  amenable.
 Hence, no F{\o}lner sequences exist on $G$ in the general case.
 However, as was shown in a recent paper \cite{Da} (see also \cite{Jo}), there is a weakly mixing infinite measure preserving $G$-action $R$ admitting an  $R$-F{\o}lner sequence of measurable subsets of finite measure. 
 Our idea is  to adapt the argument from \cite{DaKo} to the general setting of the present work by utilizing $R$ in place of the $G$-action on itself by left rotations.
 We note that in contrast with \cite{DaKo},  we do not use the ergodic property of the IDPFT actions which is based on the Schmidt-Walters theorem
 while proving Theorems~A and B for the case  $K=II_\infty$ in Section~3 (see Theorem~3.1 and Corollary~3.5). 
 Instead, we utilize an approximation technique developed in the proof of \cite{Da, Theorem~B}.

 However, this technique does not work in the other cases.
 On the other hand, the above-mentioned  theorem on ergodicity of conservative IDPFT actions 
 ---which was the crucial  tool in \cite{DaKo}---can  be applied neither because it was  proved for the case of infinite products of mildly mixing actions \cite{DaLe}, while in our case we have only weak mixing $G$-actions.
 That is why we first (in Section~1) prove a ``directional'' version of  this theorem for the weak mixing case in a similar way as it was done in \cite{ArIsMa} in the context of nonsingular Gaussian actions (see Propositions~1.16 and 1.5). 
 We also establish a ``directional'' version of the conservativeness criteria for Poisson suspensions from \cite{DaKoRo2} and \cite{DaKo} there (see Proposition~1.20). 
 Combining these results with the approach from \cite{DaKo} modified with the usage of weakly mixing actions  $R$ (see above), we prove Theorems~A and B for the case $K=III_0$ in Section~4 (see Theorems~4.1, Corollary~4.3).
 We note also that an extra argument is needed to construct $R$ whose centralizer is  ``large''---this holds automatically  in \cite{DaKo} because the $G$-action by left rotations on $G$ commutes with the  $G$-action by right rotations.

 Similar techniques combined with the argument from \cite{DaKo} are used in Section~5  to prove Theorems~A and B for the case 
 $K=III_\lambda$, $0<\lambda<1$ (see Theorem~5.1 and Corollary 5.3).
 However we also need  new additional tools  to settle this case:
 a  ``directional'' version of the famous Maharam theorem on conservativeness of the Radon-Nikodym extensions 
 (see Theorem~1.9) 
 and a theorem on  the ``symmetries'' of the $\sigma$-algebra of fixed subsets for  the Maharam extension
 (see Theorem~1.17). 
 
The remaining case $K=III_1$ of Theorems~A and B is considered in  Section~6.
Only a slight modification of the argument from Section~5 is needed to prove the theorems in this case (see Theorem~6.1 and Corollary 6.2).

\endcomment

In Section~7 we use Poisson suspensions to introduce nonsingular Bernoulli $G$-actions for arbitrary locally compact second countable group $G$.
When $G$ is discrete, these Bernoulli $G$-actions are nonsingular Bernoulli in the usual sense.
Theorem~D is proved there.
 
 In  the final Section~8 we prove Proposition~E.
 
 
 {\sl Acknowledgement.} I thank Zemer Kosloff for useful remarks and discussions and for drawing my attention to a problem raised in \cite{ArIsMa}.

\head 1. Definitions and preliminaries\endhead
\subhead Unitary representations of $G$ and Koopman representations of nonsingular actions
\endsubhead
Let $V=(V(g))_{g\in G}$ be a weakly continuous unitary representation of $G$ in a separable Hilbert space $\Cal H$. 
We will always assume that $V$ is a complexification of an orthogonal representation of $G$ in a real Hilbert space.
We recall several classical concepts from the theory of  unitary representations.

\definition{Definition 1.1} 
\roster
\item"(i)"
$V$ is called {\it weakly mixing} if $V$ has no nontrivial finite dimensional invariant subspaces.
\item"(ii)" $V$ is called {\it mixing} if $V(g)\to 0$ as $g\to\infty$ in the weak operator topology. 
\endroster
 \enddefinition

  The {\it Fock space $\Cal F(\Cal H)$ over $\Cal H$} is the orthogonal sum $\bigoplus_{n=0}^\infty\Cal H^{\odot n}$, where $\Cal H^{\odot n}$ is the $n$-th symmetric tensor power of $\Cal H$ when $n>0$ and 
 $\Cal H^{\odot 0}:=\Bbb C$.
 By $\exp V=(\exp V(g))_{g\in G}$ we denote the corresponding unitary representation of $G$ in $\Cal F(\Cal H)$, i.e.
 $\exp V(g):=\bigoplus_{n=0}^\infty V(g)^{\odot n}$ for each  $g\in G$.

 \proclaim{Fact 1.2} The following are equivalent:
  \roster
  \item"\rom{(i)}" $V$ is weakly mixing.
      \item"\rom{(ii)}" 
 There is a sequence $g_n\to\infty$ in $G$
   such that $V(g_n)\to 0$ weakly as $n\to\infty$.
    \item"\rom{(iii)}"$\exp V\restriction( \Cal F(\Cal H)\ominus\Bbb C)$ is weakly mixing.
  \endroster
 \endproclaim
 
 The equivalence (i)$\Leftrightarrow$(ii) follows from \cite{BeRo, Corollary~1.6, Theorem~1.9}.
 This equivalence imply immediately the equivalence (i)$\Leftrightarrow$(iii) (see also \cite{GlWe1, Theorem~A3}).

\subhead Nonsingular and measure preserving $G$-actions\endsubhead
Let $S=(S_g)_{g\in G}$ be a nonsingular  $G$-action  on a $\sigma$-finite standard measure space $(Z,\goth Z,\nu)$.
 Denote by $U_S=(U_S(g))_{g\in G}$ the associated (weakly continuous) unitary {\it Koopman representation} of $G$
 in $L^2(Z,\nu)$:
 $$
 U_S(g)f:=f\circ S_g^{-1}\sqrt{\frac{d\nu\circ S_g^{-1}}{d\nu}},\quad\text{for all }g\in G.
 $$
 All  actions in this paper are assumed to be effective (or faithful), i.e. $S_g=\text{Id}$ if and only
 if $g=1_G$. 
 
We recall several basic concepts related to nonsingular actions of $G$.

\definition{Definition 1.3}  
Let $H$ be a non-relatively compact subset of $G$.
\roster
\item"(i)"
$S$ is called {\it totally dissipative} if the partition of $Z$ into the $S$-orbits is measurable and the $S$-stabilizer
of a.e. point is compact, i.e. there is a measurable subset of $Z$ which meets a.e. $S$-orbit exactly once, and for a.e. $z\in Z$, the subgroup $\{g\in G\mid S_gz=z\}$ is compact in $G$.
\item"(ii)"
$S$ is called {\it conservative} if there is no any $S$-invariant subset $A\subset Z$ of positive measure such that the restriction of $S$ to $A$ is totally dissipative.
\item"(iii)"
$S$ is called {\it conservative along $H$} if for each subset $A\in\goth Z$  with $\nu(A)>0$
and each compact subset $K\subset G$,
there is $h\in H\setminus K$ such that $\nu(S_hA\cap A)>0$. 
\item"(iv)"
$S$ is called {\it ergodic} if each measurable $S$-invariant subset of $Z$ is either $\mu$-null or $\mu$-conull.
\item"(v)" $S$ is called {\it weakly mixing} if for each ergodic probability preserving 
$G$-action $R=(R_g)_{g\in G}$, the product $G$-action $(S_g\times R_g)_{g\in G}$ is 
ergodic.
\item"(vi)" $S$ is {\it of infinite ergodic index} if for each $l>0$, the the $l$-th power $(S_g^{\times l})_{g\in G}$  of $S$ is ergodic.
\item"(vii)" $S$ is called {\it of 0-type} along $H$ if $U_S(g)\to 0$ as $H\ni g\to\infty$ in the weak operator topology.
\item"(viii)" $S$ is called {\it of 0-type} if $S$ is of 0-type along $G$.
\item"(ix)" Let $S$ preserve $\nu$. A sequence $(B_n)_{n=1}^\infty$ of measurable subsets of $Z$ is called {\it $S$-F{\o}lner} if for each compact subset $K\subset G$, 
$$
\lim_{n\to\infty}\sup_{g\in K}\frac{\nu(S_gB_n\triangle B_n)}{\nu(B_n)}=0.
$$
 \item"(x)" Let $\nu(Z)=1$.
 A sequence $(B_n)_{n=1}^\infty$ of measurable subsets of $Z$ is called  {\it $S$-asymptotically invariant} if  for each compact subset $K\subset G$, 
 $$
\lim_{n\to\infty} \sup_{g\in K}\nu(B_n\triangle S_gB_n)\to 0.
 $$
     \item"(xi)" Let $\nu(Z)=1$. 
     Then $S$ is called {\it strongly ergodic} if   each $S$-asymptotically invariant sequence $(B_n)_{n=1}^\infty$ is trivial, i.e. $\lim_{n\to\infty}\nu(B_n)(1-\nu(B_n))=0$. 
 \item"(xii)" Let $\nu(Z)=1$ and $S$ preserve $\nu$.
 Then $S$ is called {\it mixing along $H$} if $\lim_{H\ni g\to\infty}\nu(S_gA\cap B)=\nu(A)\nu(B)$ for all $A,B\in\goth Z$.
\endroster
\enddefinition

   Let $ L^2_0(Z,\nu):= L^2(Z,\nu)\ominus\Bbb C=\{f\in L^2(Z,\nu)\mid\int_Zfd\nu=0\}$.
    
  \proclaim{Fact 1.4} 
     \roster
        \item"\rom{(i)}" If the $G$-action $(S_g\times S_g)_{g\in G}$ is ergodic then 
        $S$ is weakly mixing \cite{GlWe2, Theorem~1.1}.
   \item"\rom{(ii)}" Let $\nu(Z)=1$ and $S$ preserve $\nu$.
   Then $S$ is weakly mixing if and only if $U_S\restriction L^2_0(Z,\nu)$ is weakly mixing \cite{BeRo}.
 \item"\rom{(iii)}" Given a unitary representation $V$ which is the complexification of an orthogonal representation of $G$, let $(Y,\goth C,\nu, T)$ denote the probability preserving Gaussian dynamical system ($G$-action) associated with $V$.
  Then there is a canonical unitary equivalence of 
 $U_T$ and  $\exp V$ \cite{Gu}.
  \endroster
  \endproclaim

 We say that a sequence $(g_k)_{k=1}^\infty$ in $G$ is {\it dispersed} if for each compact $K\subset G$, there is $N>0$ such that if $l>m>N$ then $g_lg_m^{-1}\not\in K$.
The following statement is a part of \cite{ArIsMa, Lemma~7.15}. 
We give an alternative proof of it.

 \proclaim{Proposition~1.5} Let $S$ be a nonsingular $G$-action on a probability space $(Z,\goth Z,\nu)$.
 If $(g_k)_{k=1}^\infty$  is a dispersed sequence of $G$-elements
such that $\sum_{k=1}^\infty\frac{d\nu\circ S_{g_k}^{-1}}{d\nu}(z)=+\infty$  at a.e. $z\in Z$ then
$S$ is conservative along the subset $H:=\{g_lg_m^{-1}\mid l,m\in\Bbb N, l>m\}\subset G$.
 \endproclaim
 
 \demo{Proof} We first note that since  $(g_k)_{k=1}^\infty$  is  dispersed, $H$ is not relatively compact.
 Let $A\in\goth Z$ with $\nu(A)>0$.
 Then 
 $$
+\infty= \int_A\sum_{k=1}^\infty\frac{d\nu\circ S_{g_k}^{-1}}{d\nu}(z)\,d\nu(z)=\sum_{k=1}^\infty
\nu(S_{g_k}^{-1}A).
 $$
 Hence for each $N>0$, there exist integers $l>m>N$ such that $\nu(S_{g_l}^{-1}A\cap S_{g_m}^{-1}A)>0$.
 Therefore $\nu(A\cap S_{g_lg_m^{-1}}A)>0$.
 Since $(g_k)_{k=1}^\infty$  is a dispersed, it follows that for each compact subset $K\subset G$, there is an element $g\in H\setminus K$ such that $\nu(A\cap S_{g}A)>0$.
 Hence $S$ is conservative along $H$.
 \qed
 \enddemo

The following proposition  can be interpreted as a ``directional''  refinement of the Schmidt-Walters theorem \cite{ScWa}.
It was proved in   \cite{ArIsMa, Theorem~7.14}.
For completeness of our argument, we state the proposition here with a (modified) proof.

\proclaim{Proposition~1.6} Let $H$ be a non-relatively compact subset in $G$.
Let $R$ be a  measure preserving $G$-action on a standard probability space $(X,\goth B,\mu)$ and let $S$ be a  nonsingular $G$-action on a standard probability space $(Z,\goth Z,\nu)$.
 If $R$ is mixing along $H$ and
  $S$  is conservative along $H$ 
then each $(R_h\times S_h)_{h\in H}$-invariant subset of $X\times Z$ is the Cartesian product of $X$ with an $(R_h)_{h\in H}$-invariant subset of $Z$.
\endproclaim

\demo{Proof} Let a subset $A\subset X\times Z$ be $(R_h\times S_h)_{h\in H}$-invariant
and $\mu\otimes\nu(A)\ge0.5$.
For each $z\in Z$, we let $A_z:=\{x\in X\mid (x,z)\in A\}$ and set
 $Z_+:=\{z\in Z\mid\mu(A_z)\ge 0.5\}$.
Then
$
\kappa(Z_+)>0.
$
For each  $\epsilon>0$,
 there exist  a countable family $(B_m)_{m=1}^\infty$ of  subsets $B_m\subset Z$  and a countable family $(C_m)_{m=1}^\infty$ of  subsets $C_m\subset X$ such that 
$\nu(B_m)>0$, $\bigsqcup_{m=1}^\infty B_m=Z_+$
and $\mu(C_m\triangle A_z)<\epsilon$ for each $z\in B_m$ and all $m\in\Bbb N$.
We note that  $A$ is $(R_h\times S_h)_{h\in H}$-invariant if and only if  $A_{S_h^{-1}z}=R_hA_z$
for all $h\in H$ at a.e. $z\in Z$.
Fix $m>0$.
Select a sequence $(h_n)_{n=1}^\infty$ in $H$ such that $h_n\to\infty$ in $G$ as $n\to\infty$, $\nu(S_{h_n}B_m\cap B_m)>0$. 
Take $z_n\in S_{h_n}B_m\cap B_m$, $n\in\Bbb N$.
Then
$$
\mu(R_{h_n}C_m\cap C_m)=\mu(R_{h_n}A_{z_n}\cap A_{S_{h_n}^{-1}z_n})\pm 2\epsilon=\mu(A_{z_n})\pm 2\epsilon=\mu(C_m)\pm3\epsilon.
$$
Passing to the limit as $n\to\infty$, we obtain that 
$$
 3\epsilon\ge |\mu(C_m)^2-\mu(C_m)|\ge|\mu(C_m)-1|(0.5-\epsilon).
 $$
Hence $\mu(A_z)>1-12\epsilon$ for all $z\in B_m$.
Since $m$ and $\epsilon$ are arbitrary, we conclude that $\mu(A_z)=1$ for a.a. $z\in Z_+$.
Consider now $(X\times Z)\setminus A$ instead of $A$.
Then the similar argument yields that
$\nu(A_z)=0$ for a.a. $z\not\in Z_+$.
Hence $A=X\times Z_+$ mod~0.
Of course, $Z_+$ is $(R_h)_{h\in H}$-invariant.
\qed
\enddemo

\subhead Krieger type\endsubhead
Let $T=(T_g)_{g\in G}$ be an ergodic  nonsingular $G$-action on a standard non-atomic measure space
$(X,\goth B,\mu)$.
The {\it full group} $\big[T,\mu\big]$ consists of those $\mu$-nonsingular transformations $Q$ of $X$
for which there is a countable partition  $\Cal P$ of $X$ and a map $\Cal P\ni P\mapsto g_P\in G$ such that $Qx=T_{g_P}x$ at a.e. $x\in P$ for each $P\in\Cal P$.

If there is a $\mu$-equivalent $\sigma$-finite $T$-invariant measure then
$T$ is called {\it of type $II$}.
If the $T$-invariant measure is finite then $T$ is called {\it of type $II_1$};
if the $T$-invariant measure is infinite then $T$ is called {\it of type $II_\infty$}.
If $T$ is not of type $II$ then it is called {\it of type $III$}.
The type $III$ admits further classification into subtypes.

We first recall that {\it the Radon-Nikodym cocycle $\rho_\mu$ of $T$} is a measurable mapping 
$$
\rho_\mu:G\times X\ni (g,x)\mapsto\rho_\mu(g,x):=\frac{d\mu\circ T_g}{d\mu}(x)\in\Bbb R_+^*.
$$
An element $r$ of the multiplicative group $\Bbb R^*_+$ is called an  {\it essential value of $\rho_\mu$} if for each neighborhood $U$ of $r$ and each
subset $A\in\goth B$ of positive measure there exist a subset $B\in \goth B$ of positive measure and an element $g\in G$ such that $B\cup T_g B\subset A$ and
$\frac{d\mu\circ T_g}{d\mu}(x)\in U$ for each $x\in B$.
The set of all essential values of $\rho_\mu$ is denoted by $r(T,\mu)$.
It is a closed subgroup of $\Bbb R_+^*$.
It is easy to verify that if a measure $\gamma$ is equivalent to $\mu$ then $r(T,\mu)=r(T,\gamma)$.

If $r(T,\mu)=\Bbb R_+^*$ then $T$ is called {\it of type $III_1$}; if there is $\lambda\in(0,1)$ such that $r(T,\mu)=\{\lambda^n\mid n\in\Bbb Z\}$ then 
$T$ is called {\it of type $III_\lambda$}.
If $T$ is of type $III$ but not of type $III_\lambda$ for any $\lambda\in(0,1]$ then $T$ is called {\it of type $III_0$}.

We will  need the following folklore approximation result.

\proclaim{Fact 1.7}
Let $\goth B_0\subset \goth B$ be  a dense subring. 
Let $\delta>0$ and $s\in\Bbb R_+^*$.
If for each $A\in\goth B_0$ of positive measure  and every neighborhood $U$ of $s$ there is a subset
$B\in \goth B$  and an element $\theta\in\big[T,\mu\big]$ such that $B\cup\theta B\subset A$, $\mu(B)>\delta\mu(A)$ and 
$\frac{d\mu\circ \theta}{d\mu}(x)\in U$ for each $x\in B$ then $s\in r(T,\mu)$.
\endproclaim

\subhead Maharam extension and the associated flow\endsubhead
Let $S$ be a nonsingular $G$-action on a standard probability space $(Z,\goth Z,\nu)$.
Let $\kappa$ denote the absolutely continuous measure on $\Bbb R$ such that
$d\kappa(t)=e^{-t}dt$.
Consider the product space $(\widetilde Z,\widetilde\nu):=(Z\times\Bbb R,\nu\otimes\kappa)$.
Given $g\in G$ and $s\in\Bbb R$, we define two transformations $\widetilde S_g$ and $\widetilde s$
of $(\widetilde Z,\widetilde \nu)$
 by setting for each $(z,t)\in\widetilde Z$,
$$
{\widetilde S}_g(z,t):=\Big(S_g z, t+\log\frac{d\nu\circ S_g}{d\nu}(z)\Big)\quad\text{and}\quad
\widetilde s(z,t):=(z,t-s).
$$
Then $\widetilde S:=(S_g)_{g\in G}$ is a measure preserving $G$-action 
It is called {\it the Maharam extension} of $T$.
We note that  $(\widetilde s)_{s\in \Bbb R}$ is a totally dissipative nonsingular $\Bbb R$-action on $(\widetilde Z,\widetilde \nu)$ and
 $\widetilde S_g\widetilde s=\widetilde s\widetilde S_g$ for all $g\in G$ and $s\in\Bbb R$.
Restrict  $(\widetilde s)_{s\in\Bbb R}$ to the $\sigma$-algebra of 
$\widetilde S$-invariant subsets  of $Z\times\Bbb R$ and equip this $\sigma$-algebra with
 $\widetilde\nu$ (or, more rigorously, with a finite measure which is equivalent to $\widetilde\nu)$.
 This restriction is well defied
 as a nonsingular $\Bbb R$-action.
It is called {\it the associated flow} of $S$ and denoted by $W^S$.

\proclaim{Fact 1.8} If $(Z,\nu, S)$ is ergodic then $W^S$ is ergodic. 
Moreover,  
\roster
\item"\rom{(i)}" $S$ is of type $II$ if and only if $W^S$ is  transitive and aperiodic.
\item"\rom{(ii)}" $S$ is of type $III_\lambda$ with $0<\lambda<1$ if and only if $W^S$ is transitive but periodic with period $\log\lambda$.
\item"\rom{(iii)}" $S$ is of type $III_1$ if and only if $W^S$ is the trivial action on a singleton.
\item"\rom{(iv)}" $S$ is of type $III_0$ if and only if $W^S$ is non-transitive.
\endroster
\endproclaim

According to the  Maharam theorem,  the Maharam extension of each conservative dynamical system is conservative.
We will need the following ``directional''  refinement  of the Maharam theorem.

\proclaim{Theorem 1.9} Let  $(X,\goth B,\mu, T)$ be a nonsingular dynamical system.
Let $(g_k)_{k\in\Bbb N}$ be a dispersed sequence of $G$-elements.
If 
$$
\sum_{k=1}^\infty\Big(\frac{d\mu\circ T_{g_k}}{d\mu}(x)\Big)^{1+\alpha}=+\infty\qquad\text{ at a.e. $x$} \tag1-1
$$
for some real $\alpha\in(0,1)$
then the Maharam extension $\widetilde T$ of $T$  is conservative along the subset $\{g_lg_m^{-1}\mid l,m\in\Bbb N, l>m\}$.
\endproclaim
   \demo{Proof}
   Let $\tau$ stand for the absolutely continuous probability measure on $\Bbb R$ such that
   $d\tau(t)=\frac{\alpha} 2e^{-\alpha|t|}dt$.
   Then $\mu\otimes\tau\sim\mu\otimes\kappa$.
For each $k\in\Bbb N$ and $(x,t)\in X\times\Bbb R$, we have that
   $$
   \frac{d(\mu\otimes\tau)\circ \widetilde T_{g_k}}{d(\mu\otimes\tau)}(x,t)=
   \frac{d\mu\circ T_{g_k}}{d\mu}(x)
   e^{-\alpha\big( \big|t+\log\frac{d\mu\circ T_{g_k}}{d\mu}(x)\big|-|t|\big)}.
   $$
    We now let $N(x,t):=\Big\{k\in\Bbb N\mid \frac{d\mu\circ T_{g_k}}{d\mu}(x)\ge e^{-t}\Big\}$.
  Then we obtain that
   $$
      \frac{d(\mu\otimes\tau)\circ \widetilde T_{g_k}}{d(\mu\otimes\tau)}(x,t)=
      \cases
     e^{(|t|-t)\alpha}  \bigg(\frac{d\mu\circ T_{g_k}}{d\mu}(x)\bigg)^{1-\alpha},
     &\text{if $k\in N(x,t)$}\\
     e^{(|t|+t)\alpha}  \bigg(\frac{d\mu\circ T_{g_k}}{d\mu}(x)\bigg)^{1+\alpha},
     &\text{if if $k\not\in N(x,t)$.}
      \endcases
   $$
   If the set $N(x,t)$ is infinite then
   $$
 \sum_{k\in \Bbb N}   \frac{d(\mu\otimes\tau)\circ \widetilde T_{g_k}}{d(\mu\otimes\tau)}(x,t)\ge
  \sum_{k\in N(x,t)}     \frac{d(\mu\otimes\tau)\circ \widetilde T_{g_k}}{d(\mu\otimes\tau)}(x,t)\ge \sum_{k\in N(x,t)}  e^{-t(1-\alpha)}=\infty.
   $$
   If $N(x,t)$ is finite then \thetag{1-1} yields that
  $ \sum_{k\not\in N(x,t)}\Big(\frac{d\mu\circ T_{g_k}}{d\mu}(x)\Big)^{1+\alpha}=+\infty$ and 
hence
   $$
   \align
 \sum_{k\in \Bbb N}   \frac{d(\mu\otimes\tau)\circ \widetilde T_{g_k}}{d(\mu\otimes\tau)}(x,t)&\ge
  \sum_{k\not\in N(x,t)}     \frac{d(\mu\otimes\tau)\circ \widetilde T_{g_k}}{d(\mu\otimes\tau)}(x,t)\\
  &= e^{(|t|+t)\alpha}  \sum_{k\not\in N(x,t)}\Big(\frac{d\mu\circ T_{g_k}}{d\mu}(x)\Big)^{1+\alpha}=+\infty.
  \endalign
   $$
   Thus, in every case, i.e. for a.e. $(x,t)\in X\times\Bbb R$, we have that
   $$
   \sum_{k\in \Bbb N}   \frac{d(\mu\otimes\tau)\circ \widetilde T_{g_k}}{d(\mu\otimes\tau)}(x,t)=\infty.
   $$
   It remains to apply Proposition~1.5.
   \qed
   \enddemo

\subhead   Restricted infinite products of probability measures
\endsubhead
Let $(Y,\goth C)$ be a  standard Borel  space and let $(\gamma_n)_{n=1}^\infty$
be a sequence of probability measures on $(Y,\goth C)$.
 Let $\boldkey B:=(B_n)_{n=1}^\infty$ be a sequence of subsets from $\goth C$ such that
 $\gamma_n(B_n)>0$.
We set $(X,\goth B):=(Y,\goth C)^{\otimes\Bbb N}$.
For each $n\in\Bbb N$, let
 $
 \boldkey B^n:=Y^n\times B_{n+1}\times B_{n+2}\times\cdots\in\goth B.
 $
 Then $ \boldkey B^1\subset \boldkey B^2\subset\cdots$.
Define a measure $\gamma^{ \boldkey B}$ on  $(X,\goth B)$ by 
the following sequence of restrictions (see  \cite{Hi} for details):
 $$
 \gamma^{ \boldkey B}\restriction\boldkey B^{n}:=\frac{\gamma_1}{\gamma_1(B_1)}\otimes
 \cdots\otimes \frac{\gamma_n}{\gamma_n(B_n)}\otimes
 \frac{\gamma_{n+1}\restriction B_{n+1}}{\gamma_{n+1}(B_{n+1})}\otimes
  \frac{\gamma_{n+2}\restriction B_{n+2}}{\gamma_{n+2}(B_{n+2})}\otimes\cdots,
 $$
 $n\in\Bbb N$.
 Since these restrictions are compatible, $\gamma^{ \boldkey B}$ is well defined.
 We note that $\gamma^{ \boldkey B}$ is supported on the subset 
 $\bigcup_{n=1}^\infty \boldkey B^n\subset X$ 
 and  $\gamma^{ \boldkey B}(\boldkey B^n)=\prod_{j=1}^n\gamma_j(B_j)^{-1}$ for each $n$.
Hence, $\gamma^{ \boldkey B}$ is  $\sigma$-finite.
It is infinite if and only if $\prod_{n=1}^\infty \gamma_n(B_n)=0$. 
 
 \definition{Definition 1.10 \cite{Hi}} 
 \roster
 \item"\rom{(i)}" We call $\gamma^{ \boldkey B}$  {\it the restricted infinite product of $(\gamma_n)_{n=1}^\infty$ with respect to $\boldkey B$.}
  \item"\rom{(ii)}" 
 A  $\sigma$-finite measure $\kappa$ on  $(X,\goth B)$ is called  a {\it Moore-Hill-product (MH-product)} of $(\gamma_n)_{n=1}^\infty$ if for each $n\in\Bbb N$, there exists a $\sigma$-finite measure $\kappa_n$ on the infinite product space 
$\bigotimes_{k>n}(Y,\goth C)$ such that $\gamma=\gamma_1\otimes\cdots\otimes\gamma_n\otimes\kappa_n$.
 \endroster
 \enddefinition

Let $T=(T_g)_{g\in G}$ be a Borel action of $G$ on a Borel space $(Y,\goth C,\gamma)$.
Let $\gamma_n\circ T_g=\gamma_n$ for each $g\in G$ and each $n\in\Bbb N$.
 We define a Borel $G$-action $\boldkey T=(\boldkey T_g)_{g\in G}$ on $(X,\goth B):=(Y,\goth C)^{\otimes\Bbb N}$ by setting
$\boldkey T_g:=\bigotimes_{n=1}^\infty T_g$ for each $g\in G$.

\proclaim{Fact 1.11 \rom{(see \cite{Da, Proposition~2.6})}}
If
$\sum_{n=1}^\infty\frac{\gamma_n(B_n\triangle T_gB_n)}{\gamma_n(B_n)}<\infty$ for some $g\in G$ then
 $\boldkey T_g$ preserves $\gamma^{ \boldkey B}
 $.
 \endproclaim

 The following formula is checked straightforwardly for each $g\in G$ and $n\in\Bbb N$:
 $$
\frac{\gamma^\boldkey B(\boldkey T_g\boldkey B^n\cap\boldkey B^n)}{\gamma^{\boldkey B}(\boldkey B^n)}=\prod_{j>n}\frac{\gamma_j(T_gB_j\cap B_j)}{\gamma_j(B_j)}.\tag1-2
$$
Hence if for each compact $K\subset G$,
$$
\sum_{n=1}^\infty\sup_{g\in K}\frac{\gamma_n(B_n\triangle T_gB_n)}{\gamma_n(B_n)}<\infty
$$
then 
 the sequence $(\boldkey B^n)_{n=1}^\infty$ is $\boldkey T$-F{\o}lner.

Let $U_{j,T}$ and $U_{\boldkey T}$ denote the unitary Koopman representations of $G$ 
in $L^2(Y,\gamma_j)$ and 
$L^2(X,\gamma^{\boldkey B})$ associated with $T$ and $\boldkey T$   respectively. 

\proclaim{Fact 1.12 \rom{(cf.  \cite{Da, Lemma~2.7})}} 
Let $\sum_{n=1}^\infty\frac{\gamma_n(B_n\triangle T_gB_n)}{\gamma_n(B_n)}<\infty$.
Then for each  $n\in\Bbb N$ and two arbitrary functions $f,q\in L^2(Y^n,\bigotimes_{j=1}^n\gamma_j)$,
$$
\bigg\langle U_{\boldkey T}(g)\bigg(f\otimes \bigotimes_{j>n}1_{B_{j}}\bigg),q\otimes \bigotimes_{j>n}1_{B_{j}}\bigg\rangle
=\frac{\langle(\bigotimes_{j=1}^nU_{j,T})(g)f,q\rangle}{\prod_{j=1}^n\gamma_j(B_j)}\prod_{j>n}\frac{\gamma_j(T_gB_j\cap B_j)}{\gamma_j(B_j)}.
$$
\endproclaim

\proclaim{Corollary 1.13} Let  
$\sum_{n=1}^\infty\frac{\gamma_n(B_n\triangle T_gB_n)}{\gamma_n(B_n)}<\infty$ and
$\prod_{n=1}^\infty\gamma_n(B_n)=0$ (i.e. $\gamma^\boldkey B$ is infinite).
 If $(Y,\gamma_n,T)$ is mixing along a non-relatively compact subset  $H\subset G$ for each $n\in\Bbb N$ then  $\boldkey T$ is of 0-type along $H$.
\endproclaim
\demo{Proof} Given two positive integers $n<m$, we have that
$$
\lim_{H\ni g\to\infty}\prod_{j=n+1}^m\frac{\gamma_j(T_gB_j\cap B_j)}{\gamma_j(B_j)}=\prod_{j=n+1}^m\gamma_j(B_j)
$$
because $T$ is mixing along $H$.
Hence $\lim_{H\ni g\to\infty}\prod_{j>n}\frac{\gamma_j(T_gB_j\cap B_j)}{\gamma_j(B_j)}=0$
because $\gamma^{\boldkey B}$ is infinite.
Therefore it follows from Fact~1.12 that    $U_{\boldkey T}(g)\to 0$ weakly as $H\ni g\to\infty$. \qed
\enddemo

Given two probability measures $\alpha,\beta$ on a standard Borel space $(Y,\goth C)$,
let $\gamma$ be a third probability measure on $\goth C$ such that $\alpha\prec\gamma$ and $\beta\prec\gamma$.
The (squared) {\it Hellinger distance between $\alpha$ and $\beta$} is
$$
H^2(\alpha,\beta):=\frac12\int_Y\Bigg(\sqrt{\frac{d\alpha}{d\gamma}}-\sqrt{\frac{d\beta}{d\gamma}}\Bigg)^2=1-\int_Y\sqrt{\frac{d\alpha}{d\gamma}\frac{d\beta}{d\gamma}}\,d\gamma.
$$
This definition does not depend on the choice of $\gamma$.
The Hellinger distance is used in the Kakutani theorem on equivalence of infinite products of probability measures \cite{Ka}.
We will utilize the following fact, which is an extension of the Kakutani theorem.

\proclaim{Fact 1.14}  Let $(\gamma_n)_{n=1}^\infty$ and 
 $(\alpha_n)_{n=1}^\infty$ be  two sequences of probability measures  on $(Y,\goth C)$.
 Let $\boldkey B=(B_n)_{n=1}^\infty$ be a sequence of  subsets $B_n\in\goth C$ such that $\gamma_n(B_n)>0$. 
 Then the following are satisfied.
 \roster
 \item"(i)"  \cite{Hi, Theorem~3.6}
 $\gamma^{\boldkey B}\sim\bigotimes_{n=1}^\infty\alpha_n$ if and only if
 $\gamma_n\sim\alpha_n$ for each $n\in\Bbb N$ and
 $$\sum_{n=1}^\infty H^2\bigg(\frac 1{\gamma_n(B_n)}(\gamma_n\restriction B_n),\alpha_n\bigg)<\infty.$$
  \item"(ii)"  \cite{Hi, Theorem~3.9} Let $\gamma$ be a MH-product of $(\gamma_n)_{n=1}^\infty$.
  Then 
   there exist $a>0$ and a sequence $\boldkey D=(D_n)_{n=1}^\infty$ of subsets $D_n\in\goth C$ such that
$\gamma=a\gamma^{\boldkey D}$.
\endroster
\endproclaim

We also note  that if $\gamma^{\boldkey B}\sim\bigotimes_{n=1}^\infty\alpha_n$ then $\big(\bigotimes_{n=1}^\infty\alpha_n\big)(\boldsymbol B^k)>0$ for each $k>0$.
In particular, 
 $\prod_{n=1}^\infty\alpha_n(B_n)>0$.

\subhead  IDPFT-actions
\endsubhead
IDPFT-actions were introduced in \cite{DaLe} in the case, where $G=\Bbb Z$.
 IDPFT actions of arbitrary discrete countable groups and arbitrary locally compact Polish groups  were under consideration in  \cite{DaKo} and \cite{Da} respectively.

\definition{Definition 1.15} Let $S_n=(S_n(g))_{g\in G}$ be an ergodic  measure preserving $G$-action
on a standard probability space $(Z_n,\goth Z_n,\nu_n)$, let $\mu_n$ be a probability measure on $\goth C_n$ and let $\mu_n\sim\nu_n$ for each $n\in\Bbb N$.
We set $(Z,\goth Z,\nu):=\bigotimes_{n=1}^\infty(Z_n,\goth Z_n,\nu_n)$, $\mu:=\bigotimes_{n=1}^\infty\mu_n$,  $S(g):=\bigotimes_{n=1}^\infty S_n(g)$ for each $g\in G$ and $S:=(S(g))_{g\in G}$.
If $\mu\circ S(g)\sim\mu$ for each $g\in G$ then the nonsingular dynamical system $(Z,\goth Z,\mu, S)$
is called  an {\it infinite direct product of finite types (IDPFT)}.
\enddefinition

The Radon-Nikodym cocycle of an IDPFT system is an infinite product
 $$
 \frac{d\mu\circ S(g)}{d\mu}(z)=\prod_{n=1}^\infty\frac{d\mu_n\circ S_n(g)}{d\mu_n}(z_n)\qquad\text{ at $\mu$-a.e. $z=(z_n)_{n=1}^\infty\in Z$, $g\in G$}.
 $$

We will need the following fact, extending partly the results \cite{DaLe, Proposition~2.3}  and \cite{Da, Proposition~2.9}  about sharp weak mixing of conservative IDPFT systems with mildly mixing factors.
We extend these results to the IDPFT systems whose factors 
are  mixing along some ``directions''.

\proclaim{Proposition 1.16} Let $(Z,\goth Z,\mu, S)$ be an IDPFT system as in Definition~1.15.
Let $(g_k)_{k=1}^\infty$ be a dispersed sequence of $G$-elements such that
$\sum_{k=1}^\infty\frac{d\mu\circ S(g_k)^{-1}}{d\mu}(z)=+\infty$ at a.e. $z\in Z$.
Let $H:=\{g_lg_m^{-1}\mid l,m\in\Bbb N, l>m\}$.
If the system $(Z_n,\goth Z_n,\nu_n, S_n)$ is mixing along $H$ for each $n\in\Bbb N$  then
$(Z,\goth Z,\mu, S)$ is   weakly mixing.
\endproclaim

\demo{Proof} 
Let $R=(R_g)_{g\in G}$ be an ergodic measure preserving action of $G$ on a standard probability space $(Y,\goth Y,\kappa)$.
Then
$$
\sum_{k=1}^\infty\frac{d(\mu\otimes\kappa)\circ (S({g_k})^{-1}\times R_{g_k}^{-1})}{d(\mu\otimes\kappa)}(z,y)=\sum_{k=1}^\infty\frac{d\mu\circ S({g_k})^{-1}}{d\mu}(z)=+\infty
$$
at a.e. $(z,y)\in Z\times Y$.
It follows from Proposition~1.5 that the product $G$-action $S\times R:=(S_g\times R_g)_{g\in G}$ on the space $(Z\times Y,\mu\otimes\kappa)$
is conservative along $H$.
Let $A\subset Z\times Y$ be an $(S\times R)$-invariant subset of positive measure.
Fix $n\in\Bbb N$.
Then the product $G$-action $((\bigotimes_{j>n}S_j(g))\otimes R_g)_{g\in G}$ on
$((\prod_{j>n}Z_j)\otimes Y,(\bigotimes_{j>n}\mu_n)\otimes\kappa)$ is also conservative along $H$.
On the other hand, the $G$-action $(\bigotimes_{j=1}^nS_j(g))_{g\in G}$ on the space
$(\prod_{j=1}^nZ_j,\bigotimes_{j=1}^n\nu_n)$ is mixing along $H$.
Since the probability measure $(\bigotimes_{j=1}^n\nu_n)\otimes (\bigotimes_{j>n}\mu_n)$ on $Z$
is equivalent to $\mu$, it follows from Proposition 1.6 that there is a subset $A_n\subset (\prod_{j>n}Z_j)\times Y$ of positive measure such that $A=(\prod_{j=1}^nZ_j)\times A_n$.
Since $n$ is arbitrary, it follows that $A=Z\times B$ for some  subset $B\in\goth Y$ invariant under $R$.
Since $R$ is ergodic, $B=Y$ and hence $A=Z\times Y$.
Hence   $S\times R$ is ergodic, i.e. $S$ is weakly mixing.
\qed
\enddemo

Let Aut$(Z,\mu)$ denote the group of all invertible $\mu$-nonsingular transformations of $Z$.
We recall that {\it the weak} topology on Aut$(Z,\mu)$ is induced by the weak operator topology on the unitary group $\Cal U$ of $L^2(Z,\mu)$ via the embedding 
$$
\text{Aut}(Z,\mu)\ni R\mapsto U_R\in \Cal U.
$$ 
Then  Aut$(Z,\mu)$   is a Polish group under the weak topology.

Let $\bigoplus_{n=1}^\infty G$ stand for the direct sum of  countably many copies of $G$.
We endow this group the topology of inductive limit.
Then $\bigoplus_{n=1}^\infty G$ is $\sigma$-finite but not locally compact.
 In fact, it is not metrizable.
 Nevertheless, the nonsingular  actions of $\bigoplus_{n=1}^\infty G$ are well defined.
 Every such action $V$ is nothing but a collection $(V_n)_{n=1}^\infty$ of countable many mutually commuting
 $G$-actions defined on the  same measure space.
 We then write $V=\bigoplus_{n=1}^\infty V_n$.
Hence  the Maharam extension of $V$, we denote it by $\widetilde V$, is also well defined.
Therefore we can  construct the associated flow of $V$ in the same way as we do for the actions
of locally compact groups.

The following proposition is an analog of \cite{DaLe, Theorem~2.10} and \cite{DaKo, Proposition~1.6} for IDPFT systems with non-mildly mixing factor actions.
We provide a proof for completeness of our argument.

\proclaim{Proposition 1.17}
 Let $(Z,\goth Z,\mu, S)$ be an IDPFT system as in Definition~1.15.
Let $(g_k)_{k=1}^\infty$ be a dispersed sequence of $G$-elements such that
$$
\sum_{k=1}^\infty\Big(\frac{d\mu\circ S({g_k})^{-1}}{d\mu}(z)\Big)^\vartheta=+\infty\qquad\text{ at a.e. $z\in Z$}
$$
 for some $\vartheta\in(1,2)$.
Let $H:=\{g_lg_m^{-1}\mid l,m\in\Bbb N, l>m\}$.
If the system $(Z_n,\goth Z_n,\nu_n, S_n)$ is mixing along $H$ for each $n\in\Bbb N$
then the $\sigma$-algebra $\Cal I(\widetilde S)$ of $\widetilde S$-invariant subsets equals the $\sigma$-algebra $\Cal I\Big(\widetilde{\bigoplus_{n=1}^\infty S_n}\Big)$ of 
$\widetilde{\bigoplus_{n=1}^\infty S_n}$-invariant subsets.
Hence the associated flow of $S$ coincides with the associated flow of $\bigoplus_{n=1}^\infty S_n$.
\endproclaim
\demo{Proof} 
We first note that it follows from the condition of the proposition and Theorem~1.9 that
the Maharam extension $\widetilde S$
of $S$ is conservative along $H$.

Take a subset $A\in \Cal I(\widetilde S)$.
For every $n\in\Bbb N$, we define a measure preserving map
$E_n:(Z\times\Bbb R,\mu\otimes\eta)\to\Big(Z\times\Bbb R, \big(\bigotimes_{j=1}^n\nu_j\big)\otimes\big(\bigotimes_{j>n}\mu_n\big)\otimes\eta\Big)$ by setting
$$
E_n(z,t):=\Big(x,t+\sum_{j=1}^n\log\frac{d\mu_j}{d\nu_j}(z)\Big)
$$
It follows that for each $g\in  G$,
$$
E_n\widetilde S(g)E_n^{-1}=\bigg(\bigotimes_{j=1}^n S_j(g)\bigg)\otimes \widetilde{\bigotimes_{j>n}  S_j}(g)
$$
and that  the subset $E_nA$ is invariant under $E_n\widetilde S(g)E_n^{-1}$.
Since $\widetilde S$ is conservative along $H$,  the action $(E_n\widetilde S(g)E_n^{-1})_{g\in G}$ is also conservative along $H$.
We note that the $G$-action $\big(\widetilde{\bigotimes_{j>n}  S_j}(g)\big)_{g\in G}$
 on
 the  space $\Big(\bigotimes_{j>n}(Z_n,\mu_n)\Big)\otimes (\Bbb R,\eta)$
is
a quotient  (i.e. a factor) of
$(E_n\widetilde S(g)E_n^{-1})_{g\in G}$.
Hence $\big(\widetilde{\bigotimes_{j>n}  S_j}(g)\big)_{g\in G}$
is also conservative along $H$.
On the other hand, the measure preserving $G$-action $\Big(\bigotimes_{j=1}^n S_j(g)\Big)_{g\in G}$ on the probability space $\bigotimes_{j=1}^n(Z_n,\goth Z_n,\nu_n)$ is mixing along $H$.
Hence by Proposition~1.6, $E_nA=\Big(\bigotimes_{j=1}^n Y_j\Big)\times A_n$ for some subset $A_n\subset \Big(\bigotimes_{j>n} Y_j\Big)\times\Bbb R$.
In particular, $E_nA$ is invariant under the $G^n$-action $\Big(\bigoplus_{j=1}^n S_j\Big)\otimes I$.
Hence $A$ is invariant under $E_n^{-1}\Big(\Big(\bigoplus_{j=1}^n S_j\Big)\otimes I\Big)E_n$ which is exactly the Maharam extension of the $G^n$-action $\big(\bigoplus_{j=1}^n S_j\big)\otimes I$ on $(Z,\goth Z,\mu)$.
Since $n$ is arbitrary, $A\in {\Cal I}\Big(\widetilde{\bigoplus_{n=1}^\infty S_n}\Big)$, as desired. 

Conversely, let
$A\in {\Cal I}\Big(\widetilde{\bigoplus_{n=1}^\infty S_n}\Big)$.
We note that for each $g\in G$, the sequence of transformations  $(S_1(g)\times\cdots \times S_n(g)\times I)_{n=1}^\infty$
converges weakly to the transformation  $S(g)$ in Aut$(Z,\mu)$.
Since the mapping 
$$
\text{Aut}(Z,\mu)\ni Q\mapsto \widetilde Q\in\text{Aut}(Z\times \Bbb R,\mu\otimes\tau)
$$
 is weakly continuous, the Maharam extension of the transformation
${S(g)}$ is the weak limit of the sequence of Maharam extensions of the transformations $S_1(g)\times\cdots\times S_n(g)\times I$ as $n\to\infty$.
It follows that $A\in \Cal I(\widetilde S)$,
and we are done.
\qed
\enddemo

 \subhead Nonsingular Poisson suspension\endsubhead
 Let $(X,\goth B)$ be a standard Borel space and let $\mu$ be an infinite $\sigma$-finite 
non-atomic measure on $X$.
 Let $X^*$ be the set of purely atomic  ($\sigma$-finite) measures on $X$.
  For each subset $A\in\goth B$ with $0<\mu(A)<\infty$, we define a mapping $N_A:X^*\to\Bbb R$ by setting
 $N_A(\omega):=\omega(A)$.
 Let $\goth B^*$ stand for the smallest $\sigma$-algebra on $X^*$ such that the mappings
 $N_A$  are all $\goth B^*$-measurable.
There is a unique probability measure $\mu^*$  on $(X^*,\goth B^*)$ satisfying the following two conditions:
\roster
\item"---"  the measure $\mu^*\circ N_A^{-1}$ is the Poisson distribution
with parameter $\mu(A)$ for each $A\in\goth B$ with $\infty>\mu(A)>0$,
\item"---" given a finite family $A_1,\dots, A_q$ of mutually disjoint subsets $A_1,\dots, A_q\in \goth B$ of 
finite positive measure, the corresponding random variables $N_{A_1},\dots, N_{A_q}$ defined on the space $(X^*,\goth B^*, \mu^*)$ are independent.
\endroster
 Then $(X^*,\goth B^*, \mu^*)$ is a Lebesgue space.
 Given  a subset $B\in\goth B$ and an integer $n\in\Bbb Z_+$, we denote by
$[B]_n$ the cylinder $\{\omega\in X^*\mid \omega(B)=n\}$.
 We now
 let
$$
\text{Aut}_1(X,\mu):=\bigg\{S\in \text{Aut}(X,\mu)\,\Big|\, {\frac{d\mu\circ S}{d\mu}}-1\in L^1(X,\mu)\bigg\}.
$$ 
If $S\in \text{Aut}_1(X,\mu)$, we put $\chi(S):=\int_X( {\frac{d\mu\circ S}{d\mu}}-1)d\mu$.
Then Aut$_1(X,\mu)$ is  a subgroup of Aut$(X,\mu)$ and $\chi$ is a homomorphism
of Aut$_1(X,\mu)$ onto $\Bbb R$.
Suppose that  $T=(T_g)_{g\in G}$ is a nonsingular  $G$-action on $(X,\goth B,\mu)$. 
We now define a Borel transformation $T_g^*$ of $X^*$ by setting
$$
T_g^*\omega:=\omega\circ T_g^{-1} \qquad\text{for all $\omega\in X^*$ for each $g\in G$.}
$$

\proclaim{Fact 1.18}
\roster
\item"\rom{(i)}" 
If $T_g\in\text{\rom{Aut}}_1(X,\mu)$ for each $g\in G$ then
 $T^*:=(T_g^*)_{g\in G}$  is a well defined nonsingular $G$-action on $(X^*,\goth B^*, \mu^*)$ \cite{DaKoRo1, \S 6.1 and \S4}.
\item"\rom{(ii)}"
If $\nu$ is a $\sigma$-finite measure on $(X,\goth B)$ then $\nu^*\sim\mu^*$ if and only if
$\mu\sim\nu$ and $\sqrt{\frac{d\mu}{d\nu}}-1\in L^2(X,\nu)$ (see \cite{Ta} or \cite{DaKoRo1, Theorem~3.3}).
\item"\rom{(iii)}"
  Moreover,  if  $\frac{d\mu}{d\nu}-1\in L^1(X,\nu)$ then
$$
\frac{d\mu^*}{d\nu^*}(\omega)=e^{-\int_X(\frac{d\mu}{d\nu}-1)d\nu}\prod_{\omega(\{x\})=1}\frac{d\mu}{d\nu}(x)\qquad\text{ at a.e. $\omega\in X^*$.}
$$
\item"\rom{(iv)}"
 It follows (see also \cite{DaKoRo1, Corollary~4.1(3)}) that  if $S\in\text{\rom{Aut}}_1(X,\mu)$ and $\chi(S)=0$ then 
  $$
  \frac{d\mu^*\circ (S^*)^{-1}}{d\mu^*}(\omega)=\prod_{\{x\in X\mid\omega(\{x\})=1\}}\frac{d\mu\circ S^{-1}}{d\mu}(x)\qquad\text{for a.e. $\omega\in X^*$}.
  $$
  \item"\rom{(v)}" If $X$ is partitioned into countably many $T$-invariant subsets $X_n$, $n\in\Bbb N$, then  $(X^*,\mu^*,T^*)$ is canonically isomorphic to the direct product
  $\bigotimes_{n=1}^\infty(X_n^*,\mu_n^*,T_n^*)$, where $\mu_n$ and $T_n$ denote the restriction of $\mu$ and $T$  respectively to $X_n$.
  \endroster
  \endproclaim
  
  In particular, if $T$ preserves $\mu$ then $T^*$ preserves $\mu^*$.

\definition{Definition 1.19}
The dynamical     system $(X^*,\goth B^*, \mu^*, T^*)$ is called
 {\it the nonsingular Poisson suspension} of $(X,\goth B, \mu, T)$.\footnote{We assume that the condition of Fact~1.18(i) holds.}
 A  nonsingular $G$-action
  is called {\it Poisson} if it  is isomorphic to the Poisson suspension of some nonsingular $G$-action (see \cite{DaKoRo1} for details).
  \enddefinition

  The following proposition is an adaptation of  \cite{DaKo, Lemma~1.3} (see also \cite{DaKoRo2, Theorem~3.4}) to the case of locally compact group actions.

\proclaim{Proposition 1.20}
 Let
 $T=(T_g)_{g\in G}$ be a nonsingular $G$-action on a $\sigma$-finite measure standard non-atomic measure space $(X,\goth B,\mu)$.
 Suppose that  $T_g\in \text{\rom{Aut}}_1(X,\mu)$ and  $\chi(T_g)=0$  and
 for all $g\in G$.
 Let $(g_k)_{k\in\Bbb N}$ be a dispersed sequence of $G$-elements  such that
 $\Big(\frac{d\mu}{d\mu\circ T_{g_k}^{-1}}\Big)^2-1\in L^1(X,\mu)$.
 If there is a sequence $(b_k)_{k=1}^\infty$ of positive reals such that $b_k\le 1$ for each $k$, 
 $\sum_{k=1}^\infty b_k=+\infty$ but 
 $$
 \sum_{k=1}^\infty b_k^{1+\vartheta}e^{\int_X\Big(\Big(\frac{d\mu}{d\mu\circ T_{g_k}^{-1}}\Big)^2-1\Big)d\mu}<+\infty
 \tag1-3
 $$
 for some real $\vartheta>0$
then
$$
\sum_{k=1}^\infty\bigg(\frac{d\mu^*\circ (T_{g_k}^*)^{-1}}{d\mu^*}(\omega)\bigg)^{\frac 2{1+\vartheta}}=+\infty\qquad\text{at a.e. $\omega$.}
\tag1-4
$$
If $\vartheta\le 1$ then
 the Poisson suspension $T^*:=(T_g^*)_{g\in G}$ of $T$ is conservative along
the subset $\{g_lg_m^{-1}\mid l,m\in\Bbb N, l>m\}$.
 \endproclaim

\demo{Proof} As in the the proof of  \cite{DaKo, Lemma~1.3}, 
it follows from the 
 assumptions of the lemma that for each $k\in\Bbb N$,
$$
M_k :=
\bigg\|\frac{d\mu^*}{d\mu^*\circ (T_{g_k}^*)^{-1}}\bigg\|_2^2 =
e^{\int_X\Big(\Big( \frac{d\mu}{d\mu\circ T_{g_k}^{-1}}\Big)^2-1\Big)d\mu}.
$$
By Markov's inequality, 
$$
\align
\mu^*\bigg(\bigg\{\omega\in X^*\,\bigg|\, \frac{d\mu^*}{d\mu^*\circ (T_{g_k}^*)^{-1}}(\omega)>b_k^{-\frac {1+\vartheta}2}
\bigg\}\bigg)
&\le b_k^{1+\vartheta} M_k\\
&=b_k^{1+\vartheta} e^{\int_X\Big(\Big( \frac{d\mu}{d\mu\circ T_{g_k}^{-1}}\Big)^2-1\Big)d\mu}.
\endalign
$$
It now follows from \thetag{1-3} and the Borel-Cantelli lemma  that
$$
\frac{d\mu^*}{d\mu^*\circ (T_{g_k}^*)^{-1}}(\omega)\le b_k^{-\frac {1+\vartheta}2}
\tag1-5
$$
 for all but finitely many $k\in\Bbb N$ at a.e. $\omega$.
We can rewrite \thetag{1-5} as 
$$
\bigg(\frac{d\mu^*\circ (T_{g_k}^*)^{-1}}{d\mu^*}(\omega)\bigg)^{\frac{2}{1+\vartheta}}\ge b_k.
$$
As $\sum_{k=1}^\infty b_k=\infty$, \thetag{1-4} follows.

If  $\vartheta\le1$ then  $\frac{1+\vartheta}2\le 1$ and we deduce from \thetag{1-4} that
$$
\sum_{k=1}^\infty\frac{d\mu^*\circ (T_{g_k}^*)^{-1}}{d\mu^*}(\omega)
\ge \sum_{k=1}^\infty b_k^{\frac{1+\vartheta}2}\ge
 \sum_{k=1}^\infty b_k
=+\infty
$$
at a.e. $\omega$.
It remains to apply Proposition~1.5.
 \qed
 \enddemo

Denote by $U_T$ and $U_{T^*}$ the corresponding unitary Koopman representations of $G$ 
in $L^2(X,\mu)$ and $L^2(X^*,\mu^*)$ respectively.

\proclaim{Fact 1.21 \rom{(see \cite{Ro})}} There is a canonical unitary isomorphism of $U_{T^*}$ and $\exp U_T$.
\endproclaim

\subhead{Amenability of groups and actions}\endsubhead
Let $T=(T_g)_{g\in G}$ be a nonsingular $G$-action on standard probability space $(X,\goth B,\mu)$.
Denote by $\lambda_G$  a left Haar measure on $G$.

\definition{Definition 1.22} 
\roster
\item"(i)" $T$ is called {\it amenable} (\cite{Zi}, \cite{AD}) if there is
a $G$-invariant mean (i.e. a positive unital linear map) from $L^\infty(X\times G,\mu\otimes\lambda_G)$ onto $L^\infty(X,\mu)$, where $X\times G$ is endowed with the diagonal $G$-action (with $G$ acting on $G$ by left rotations).
\item"(ii)" $T$ is called {\it amenable in the Greenleaf sense (\cite{Gr}}, \cite{AD}) if there is a $G$-invariant mean on $L^\infty(X,\mu)$.
\endroster
\enddefinition
Of course, each probability preserving $G$-action is amenable in the Greenleaf sense. 
We also recall  the definition of  weak containment for unitary representations.

\definition{Definition 1.23} Let $\pi_1$ and $\pi_2$ be 
 two unitary representations of $G$ in Hilbert spaces $\Cal H_1$ and $\Cal H_2$ respectively.   {\it $\pi_1$ is weakly contained in $\pi_2$} if  given a vector $h\in\Cal H_1$,  a compact subset $K\subset G$ and $\epsilon>0$, there exist  finitely many vectors
$h_1,\dots, h_k\in\Cal H_2$ such that
$$
\sup_{g\in K}| \langle \pi_1(g)h,h\rangle-\sum_{j=1}^k
 \langle \pi_2(g)h_j,h_j\rangle|<\epsilon.
$$
\enddefinition

Let $L_G$ stand for the left regular representation of $G$.

\proclaim{Fact 1.24} Let $T$ be a nonsingular $G$-action.
\roster
\item"\rom{(i)}" If $G$ is amenable then  $T$ is amenable \cite{Zi}.
\item"\rom{(ii)}" If $T$ is amenable then  $U_T$ is weakly contained in $L_G$ \cite{AD, Theorem~4.3.1}.
\item"\rom{(iii)}" $U_T$ has almost invariant vectors if and only if the trivial representation of $G$ is weakly contained in $U_T$.
\item"\rom{(iv)}" $G$ is amenable if and only if  the trivial representation of $G$ is weakly contained in $L_G$.
\item"\rom{(v)}" $T$ is amenable in the Greenleaf sense if and only if 
the trivial representation of $G$ is weakly contained in $U_T$ \cite{AD, Proposition~4.1.1}.
\endroster
\endproclaim

The converse to (i) does not hold. 
For each locally compact second countable group $G$, the action of $G$ on itself by left translations is amenable.
The converse to~(ii) does not hold either \cite{AD}.

\proclaim{Proposition 1.25} Let  $(Z,\goth Z,\mu,S)$ be an IDPFT system as in Definition~1.15.
Then $S$ is amenable in the Greenleaf sense. 
Hence $S$ is amenable if and only if $G$ is amenable.
\endproclaim
\demo{Proof}
We will use the notation from Definition~1.15.
Let $\xi_n:=\frac{d\mu_n}{d\nu_n}$ for each $n\in\Bbb N$.
We will  show that the Koopman representation $U_S$ has almost invariant vectors.
Fix a compact subset $K\subset G$ and a real $\epsilon>0$.
Since $\min_{g\in K}\langle U_S(g)1,1\rangle>0$ and $\langle U_S(g)1,1\rangle=
\prod_{k=1}^\infty\langle U_{S_k}(g)1,1\rangle$, there is  $N\in\Bbb N$ such that
$$
\sup_{g\in K}\bigg|\prod_{k>N}\langle U_{S_k}(g)1,1\rangle-1\bigg|<\epsilon.
$$
Let $v_n:=\frac 1{
\sqrt{\xi_n}}$.
It is straightforward  to verify $v_n\in L^2(Z_n,\mu_n)$,  $\|v_n\|_2=1$ and 
$
\langle U_{S_n}(g)v_n,v_n\rangle=1
$
for each $g\in G$.
Let $w_n:=v_1\otimes\cdots\otimes v_n\otimes 1\otimes 1\otimes\cdots$.
Then $w_n\in L^2(Z,\mu)$ and $\|w_n\|_2=1$ for each $n\in\Bbb N$.
Moreover,
$$
\langle U_S(g)w_N,w_N\rangle=\prod_{k>N}\langle U_{S_k}(g)1,1\rangle=1\pm\epsilon
$$
for each $g\in K$.
Hence $\|U_S(g)w_N-w_N\|^2<2\epsilon$.
Thus, $U_S$ has  almost invariant vectors.
In view of Fact~1.24 (iii) and (v), $S$ is amenable in the Greenleaf sense.

 If $G$ is amenable then $S$ is amenable by Fact 1.24(i).

Suppose now that $G$ is non-amenable but $S$ is amenable.
Since $U_S$ has almost invariant vectors, it follows from 
 Fact~1.24(iii) that the trivial representation of $G$ is weakly contained in $U_S$.
Since $S$ is amenable, $U_S$ is weakly contained in $L_G$ by Fact~1.24(ii).
Hence the trivial representation of $G$ is weakly contained in $L_G$, i.e. $G$ is amenable by Fact~1.24(iv), a contradiction.
\qed
\enddemo

\head 2. General construction
\endhead

Let $S=(S_g)_{g\in G}$ be a measure preserving $G$-action on an infinite $\sigma$-finite standard measure space
$(Y,\goth Y,\kappa)$ and let $\boldkey F:=(F_n)_{n=1}^\infty$ be an $S$-F{\o}lner sequence of measurable subsets in $Y$.
Let  $(a_n)_{n=1}^\infty$ and $(\lambda_n)_{n=1}^\infty$ be two  sequences of positive reals.
For each $n\in\Bbb N$, we set 
$$
f_n:=1+(\lambda_n-1)\cdot 1_{F_n}.
$$
Define two  measures $\nu_n$ and $\mu_n$ on $Y$ by setting: $\nu_n:=\frac{a_n}{\kappa(F_n)}\,\kappa$ and $\mu_n\sim\nu_n$ with $\frac{d\mu_n}{d\nu_n}:=f_n$.
Then $\nu_n(F_n)=a_n$ and $\mu_n(F_n)=\lambda_na_n$.
We now let 
$$
(X,\goth B,\nu):=\bigsqcup_{n=1}^\infty(Y,\goth Y,\nu_n)\qquad\text{ and}\qquad  \mu:=\bigsqcup_{n=1}^\infty\mu_n.
$$
Denote by $T=(T_g)_{g\in G}$ the $G$-action on $X$ whose restriction to every copy of $Y$ in $X$ is $S$.
Then $\nu$ and $\mu$ are $\sigma$-finite measures on $(X,\goth B)$, $\nu\circ T_g=\nu$ and $\mu\circ T_g\sim\mu$ for each $g\in G$.

\proclaim{Lemma 2.1}
$T_g\in\text{\rom{Aut}}_1(X,\mu)$ if and only if
$$
\sum_{n\in\Bbb N}a_n|\lambda_n-1|\frac{\kappa
(S_gF_n\triangle F_n)}{\kappa(F_n)}<\infty.\tag2-1
$$
Moreover, $\chi(T_g)=0$ whenever
$T_g\in\text{\rom{Aut}}_1(X,\mu)$.
\endproclaim

\demo{Proof} We compute that
$$
\align
\int_X\bigg|\frac{d\mu\circ T_g^{-1}}{d\mu}-1\bigg|d\mu&=
\sum_{n\in\Bbb N}\int_Y|f_n\circ S_g^{-1}- f_n|\frac {a_n}{\kappa(F_n)}d\kappa\\
&=
\sum_{n\in\Bbb N}\int_Y|\lambda_n-1|\cdot |1_{F_n}\circ S_g^{-1}-1_{F_n}|\frac{a_n}{\kappa(F_n)}d\kappa\\
&=
\sum_{n\in\Bbb N}a_n|\lambda_n-1|\frac{\kappa
(S_gF_n\triangle F_n)}{\kappa(F_n)}.
\endalign
$$
Hence \thetag{2-1} follows.
In a similar way, we obtain that that 
$$
\chi(T_g)=\int_X\bigg(\frac{d\mu\circ T_g^{-1}}{d\mu}-1\bigg)d\mu=
\sum_{n\in\Bbb N}\frac{a_n}{\nu(F_n)}\int_Y(f_n\circ S_g^{-1}- f_n)d\nu=0.
$$
\qed
\enddemo

For each $\lambda>0$, we let $c(\lambda) :=\lambda^3-\lambda+\lambda^{-2}-1 = (1-\lambda^2)(1-\lambda^3)\lambda^{-2}\ge 0.$ 

\proclaim{Lemma 2.2}
 $\Big(\frac{d\mu}{d\mu\circ T_g^{-1}}\Big)^2-1\in L^1(X,\mu)$ if and only if
$$
\sum_{n=1}^\infty
a_n \frac{|1-\lambda_n^2|(1+\lambda_n^3)}{\lambda_n^2}\frac{\kappa
(S_gF_n\triangle F_n) }{\kappa(F_n)}
<\infty
$$
and
$$
\int_X\Big(\Big(\frac{d\mu}{d\mu\circ T_g^{-1}}\Big)^2-1\Big)d\mu=0.5\sum_{n=1}^\infty a_nc(\lambda_n)\frac{\kappa
(S_gF_n\triangle F_n) }{\kappa(F_n)}.
$$
\endproclaim

\demo{Proof}
We first observe that
$$
\align
\int_X\bigg|\bigg(\frac{d\mu}{d\mu\circ T_g^{-1}}\bigg)^2-1\bigg|d\mu
&=\sum_{n=1}^\infty
\frac{a_n}{\kappa(F_n)}\int_Y\bigg|\frac{f^3}{(f\circ T_g^{-1})^2}-f\bigg|d\kappa\\
&=\sum_{n=1}^\infty
a_n \frac{|1-\lambda_n^2|(1+\lambda_n^3)}{\lambda_n^2}\frac{\kappa
(S_gF_n\triangle F_n) }{\kappa(F_n)}
\endalign
$$
and the first claim of Lemma~2.2 follows.
In a similar way,
$$
\align
\int_X\bigg(\bigg(\frac {d\mu}{d\mu\circ T_g^{-1}}\bigg)^2-1\bigg) d\mu
&=\sum_{n=1}^\infty \frac{a_n}{\kappa(F_n)}\Big(
(\lambda_n^3-\lambda_n)\kappa(F_n\setminus S_gF_n)\\
&+(\lambda_n^{-2}-1)\kappa(S_g F_n\setminus F_n)\Big)\\
&=\sum_{k=1}^\infty a_nc(\lambda_n)\frac{\kappa
(S_gF_n\triangle F_n) }{2\kappa(F_n)},
\endalign
$$
as desired.
\qed
\enddemo

Let $(Y^*,\goth Y^*,\mu_n^*,S^*)$ denote the  Poisson suspension of the nonsingular dynamical system $(Y,\goth Y,\mu_n,S)$, $n\in\Bbb N$.
Suppose that \thetag{2-1} holds.
Then the Poisson suspension $T^*=(T_g^*)_{g\in G}$ of $T$ is well defined as a nonsingular $G$-action on a standard probability space $(X^*,\goth B^*,\mu^*)$.
By Fact~1.18(v), we have that
$(X^*,\mu^*,T^*)$ is canonically isomorphic to the direct product $\bigotimes_{n\in\Bbb N}(Y^*,\mu_n^*,S^*)$.
Since $\frac{d\mu_n}{d\nu_n}-1=f_n-1=(\lambda_n-1)1_{F_n}\in L^1(Y,\nu_n)$, 
we obtain that   $\mu_n^*\sim\nu_n^*$.
Moreover $\nu_n^*$ is invariant under $S^*$.
Thus, we have shown that  $(X^*,\goth B^*,\mu^*, T^*)$ is IDPFT.

\proclaim{Proposition 2.3 \rom{(cf.  \cite{Da, Corollary C})}}
If there are $c>0$ and an $S$-F{\o}lner sequence $(B_n)_{n=1}^\infty$  in $(Y,\goth Y,\kappa)$ such that $\kappa(B_n)=c$ for each $n\in\Bbb N$ then 
 $T^*$ is not strongly ergodic.
\endproclaim
\demo{Proof} Since
$(X^*,\mu^*,T^*)$ is isomorphic to $\bigotimes_{n\in\Bbb N}(Y^*,\mu_n^*,S^*)$ and $\mu_n^*\sim\nu_n^*$ for each $n\in\Bbb N$,
it follows that 
$(Y^*,\nu_1^*,S^*)$ is a factor of
$(X^*,\mu^*,T^*)$.
Hence it suffices to show that 
$(Y^*,\nu_1^*,S^*)$ is not strongly ergodic.
We note that
 $\nu_1^*([B_n]_0)=e^{-\nu_1(B_n)}$ and
$$
\nu_1^*([B_n]_0\cap S_g^*[B_n]_0)=\nu_1^*([B_n\cup S_gB_n]_0)=e^{-\nu_1(B_n\cup S_gB_n)}
$$
for each $n\in\Bbb N$.
Since $(B_n)_{n=1}^\infty$ is $S$-F{\o}lner,
 $\sup_{g\in K}|\frac{\nu_1(B_n\cup S_gB_n)}{\nu_1(B_n)}-1|\to 0$ as $n\to\infty$ for each compact subset $K\subset G$.
 Hence 
 $\sup_{g\in K}|\nu_1(B_n\cup S_gB_n)-\nu_1(B_n)|\to 0$ as $n\to\infty$.
It follows that $([B_n]_0)_{n=1}^\infty$ is an  asymptotically invariant  sequence in $(Y^*,\goth Y^*,\nu_1^*)$.
Since $\nu_1(B_n)=\frac{a_1c}{\kappa(F_1)}$, the sequence $([B_n]_0)_{n=1}^\infty$ is
nontrivial.
Hence $T^*$ is not strongly ergodic.
\qed
\enddemo

\proclaim{Proposition 2.4} Let $S$ be effective, i.e. $S_g\ne\text{\rom{Id}}$ for each $g\in G$.
If there is a measure preserving transformation $Q$ of $(Y,\goth Y,\kappa)$ such that
\roster
\item"\rom{(i)}"
$QS_g=S_gQ$ for each $g\in G$ and
\item"\rom{(ii)}"
if $QA=A$ for a subset $A\in\goth Y$ then either $\kappa(A)=0$ or $\kappa(A)=+\infty$
\endroster
then $T^*$ is free.
\endproclaim
\demo{Proof}
Since $(X^*,\mu^*,T^*)=\bigotimes_{n=1}^\infty (Y^*,\mu_n^*,S^*)$
and  $\mu_1^*\sim\nu_1^*$ by Fact~1.18(ii),
it suffices to show that the action
$S^*$ on $(Y^*,\nu_1^*)$ is free.
Since $Q$ preserves $\nu_1$, it follows that the Poisson suspension $Q^*$ of $Q$
preserves $\nu_1^*$.
We deduce from (ii) that $Q^*$ is ergodic. 
Moreover, (i) yields that $Q^*$ commutes with $S_g^*$ for each $g\in G$.

 Denote by $\Cal G$ the space of all closed subgroups of $G$ and endow $\Cal G$ with the Fell topology \cite{Fe}.
 Then $\Cal G$ is a compact metric space \cite{Fe}.
 Given $\omega\in Y^*$, let $G_\omega:=\{g\in G\mid S^*_g\omega=\omega\}$
 stand for  the stability group of $S^*$ at $\omega$.
 Then $G_\omega\in\Cal G$ at each $\omega$ \cite{AuMo, I, Proposition~3.7}.
  The mapping
 $
 \eta:Y^*\ni\omega\mapsto G_\omega\in \Cal G
 $
is  Borel \cite{AuMo, II, Proposition~2.3}.
Since $Q^*$ commutes with $S^*$, it is straightforward to verify that  $\eta$ is invariant under $Q^*$.
As $Q^*$ is ergodic,  we obtain that $\eta$ is constant.
Thus,  there is a subgroup $H\in\Cal G$ such that $G_\omega=H$ at a.e. $\omega\in Y^*$.
Therefore $S^*_g=\text{Id}$ for each $g\in H$.
This implies that $S_g=\text{Id}$ for each $g\in H$.
Since $S$ is effective, $H$ is trivial.
Thus, we obtain that   $S^*$  (and hence $T^*$) is free, as desired.
\qed
\enddemo

\head 3. Type $II_\infty$ nonsingular  Poisson  suspensions
\endhead
In this section we prove the implications (1)$\Rightarrow$(3) of Theorems~A and B for  $K=II_\infty$.
Thus, we assume that  $G$ does not have property (T).
 Then there exists a measure preserving  $G$-action
$S=(S_g)_{g\in G}$ on     an infinite $\sigma$-finite standard measure space $(Y,\goth Y,\kappa)$ such that
\roster
\item"$(\alpha_1)$"
there is an $S$-F{\o}lner sequence $\boldkey F=(F_n)_{n=1}^\infty$ such that $\kappa(F_n)=1$ for each $n\in\Bbb N$,
\item"$(\alpha_2)$" the corresponding unitary Koopman representation $U_S$ of $G$ is weakly mixing
\endroster 
(see \cite{Jo} or \cite{Da, Theorem~D(ii)}).
We can assume without loss of generality that $S$ is free.
Indeed, let $L=(L_g)_{g\in G}$ stand for the $G$-action on $G$ by left  translations. 
We endow $G$ with the left Haar measure.
Then $L$ is an infinite measure preserving free $G$-action.
The Poisson suspension $L^*=(L^*_g)_{g\in G}$ of $L$ is a free mixing probability preserving 
$G$-action \cite{OrWe}.
If $S$ is not free then we replace $S$ with  the product $G$-action $(S_g\times L^*_g)_{g\in G}$
which is free and satisfies  $(\alpha_1)$ and $(\alpha_2)$. 

Passing to a subsequence in $\boldkey F$, we may (and will) also assume without loss of generality that
\roster
\item"$(\alpha_3)$" $\sum_{n=1}^\infty\sup_{g\in K}{\kappa(S_gF_n\triangle F_n)}
<+\infty$ for each compact $K\subset G$.
\endroster
Fix two sequences $(a_n)_{n=1}^\infty$ and $(\lambda_n)_{n=1}^\infty$
of positive reals such that
\roster
\item"$(\alpha_4)$" $a_n:=1$ for all $n\in\Bbb N$ and  $\sum_{n=1}^\infty\lambda_n<+\infty$.
\endroster
We will also need one extra property:
\roster
\item"$(\alpha_5)$" There is a $\kappa$-preserving transformation $Q$  that commutes with $S_g$ for each $g\in G$ and such that each $Q$-invariant subset is either of 0-measure or of infinite measure.
\endroster

\remark{Remark \rom{3.1}}
If there exists an action $S$ satisfying $(\alpha_1)$--$(\alpha_4)$, it is easy to construct a new $G$-action $S'$ satisfying  $(\alpha_1)$--$(\alpha_5)$.
Indeed, just put $Y':=Y\times \Bbb Z$, $\kappa'=\kappa\otimes\delta_{\Bbb Z}$, 
$F_n':=F_n\times\{0\}$,
$S'_g:=S_g\otimes\text{Id}$ and $Q=\text{Id}\otimes r$, where $\delta_{\Bbb Z}$ is the counting measure on $\Bbb Z$ and $r$ is the translation by $1$ on $\Bbb Z$.
It is straightforward to verify that  $(Y',\kappa', S', (F_n')_{n=1}^\infty, Q)$
satisfy  $(\alpha_1)$--$(\alpha_5)$.
\endremark

Denote by $(X,\goth B,\mu, T)$ the dynamical system associated with $(Y,\goth Y,\kappa,S)$,  $\boldkey F$, $(a_n)_{n=1}^\infty$ and $(\lambda_n)_{n=1}^\infty$ via the general construction described in  \S 2.
We deduce from  Lemma~2.1, $(\alpha_3)$ and  $(\alpha_4)$ that $T_g\in\text{Aut}_1(X,\mu)$ for each $g\in G$.
Hence the nonsingular Poisson suspension $(X^*,\goth B^*,\mu^*,T^*)$ of $(X,\goth B,\mu, T)$ is well defined according to Fact~1.18(i).

\proclaim{Theorem 3.2} 
Let $S=(S_g)_{g\in G}$ be as above.
Then there is a sequence $k_n\to\infty$ such that the   general construction of \S2 applied to $(Y,\goth Y,\kappa,S)$,   $(F_{k_n})_{n=1}^\infty$, $(a_{k_n})_{n=1}^\infty$ and
 $(\lambda_{k_n})_{n=1}^\infty$ yields the dynamical system 
 $(\widetilde X,\widetilde{\goth B},\widetilde\mu, \widetilde T)$ whose 
  Poisson suspension $(\widetilde X^*,\widetilde{\goth B}^*, \widetilde\mu^*, \widetilde T^*)$  is  free, of infinite ergodic index and
   of Krieger type $II_\infty$, non-strongly ergodic and  IDPFT.
 Hence, $ \widetilde T^*$ is  amenable in the Greenleaf sense.
Moreover, $ \widetilde T^*$ is amenable if and only if $G$ is amenable.
\endproclaim

Before we get  to the proof  of this  theorem, let us first state  an  approximation lemma.

 \proclaim{Lemma 3.3 \cite{Da, Proposition 3.1}} Let $R=(R_g)_{g\in G}$ be a measure preserving $G$-action  on an infinite   $\sigma$-finite standard measure space $(Z,\goth Z,\tau)$.
  Let $\goth Z_0\subset \goth Z$ stand for the ring of subsets of finite measure.
Let  $Z_1\subset Z_2\subset\cdots$ be a sequence of subsets from $\goth Z_0$ such that
  for each $n>0$, there is a finite partition $\Cal P_n$ of $Z_n$
into subsets of equal measure
satisfying the following:
\roster
\item"---" $\bigcup_{n=1}^\infty Z_n=Z$ and
\item"---"
$(\Cal P_n)_{n=1}^\infty$ approximates $(\goth Z_0,\tau)$, i.e.
for each  $B\in \goth Z_0$ and $\epsilon>0$, there is $N>0$ such that if $n>N$ then there is  a $\Cal P_n$-measurable subset $B_n$ with $\tau(B\triangle B_n)<\epsilon$.
\endroster
 If   for each $n>0$ and every pair of $\Cal P_n$-atoms $A$ and $B$,
  there exist a finite family
  $g_1,\dots,g_l\in G$ and  mutually disjoint subsets\footnote{Some of these subsets can be of $0$ measure.} $A_1,\dots,A_l$ of $A$ such that the subsets $R_{g_1}A_1,\dots,R_{g_l}A_l$ are mutually disjoint, 
  $\bigsqcup_{i=1}^lR_{g_i}A_i\subset B$    and $\tau(\bigsqcup_{i=1}A_i)>0.1\tau(A)$ then
    $R$ is ergodic.
    \endproclaim

\demo{Proof of Theorem~3.2}
As we noted in \S2, the dynamical system 
$(X^*,\goth B^*,\mu^*,T^*)$ is canonically isomorphic to the infinite direct product  $\bigotimes_{n\in\Bbb N}(Y^*,\goth Y^*,\mu_n^*,S^*)$.
Let $B_n:=[F_n]_0$ and $\boldkey B:=( B_n)_{n=1}^\infty$.

{\sl Claim A.}
$
\sum_{n=1}^\infty H^2\Big(\frac1{\nu_n^*(B_n)}\nu_n^*\restriction B_n,\mu_n^*\Big)<\infty.
$

We first recall that $\nu_n:=\frac{a_n}{\kappa(F_n)}\,\kappa=\kappa$ in view of 
$(\alpha_1)$ and  $(\alpha_4)$.
Therefore
$\nu_n(F_n)=1$ and  $\nu_n^*(B_n)=e^{-\nu_n(F_n)}=e^{-1}$.
It follows now from~Fact~1.18(iii) that if $\omega\in B_n$ then
$$
\frac{d\mu_n^*}{d\nu_n^*}(\omega)=e^{-(\lambda_n-1)\nu_n(F_n)}\sum_{k=0}^\infty\lambda_n^k1_{[F_n]_k}(\omega)=e^{1-\lambda_n}.
$$
Hence
$$
\align
H^2\bigg(\frac1{\nu_n^*(B_n)}\nu_n^*\restriction B_n,\mu_n^*\bigg)&
=1-\int_{B_n} \frac{1}{\sqrt{\nu_n^*(B_n)}}\sqrt{\frac{d\mu^*_n}{d\nu_n^*}} \,d\nu_n^*
\\
&=1- e^{-\frac12}e^{\frac{1-\lambda_n}2} \\
&=1-e^{-\frac{\lambda_n}{2}}.
\endalign
$$
Since $\sum_{n=1}^\infty\lambda_n<\infty$ by $(\alpha_4)$, Claim~A follows.

 Fact~1.14(i) and Claim~A yield that $\mu^*$
is equivalent to the restricted infinite product of $(\nu_n^*)_{n=1}^\infty$ with respect to $\boldkey B$.
We will denote this restricted product by $(\nu^*)^{\boldkey B}$.
Since $\prod_{n=1}^\infty\nu_n^*(B_n)=0$, it follows that $(\nu^*)^{\boldkey B}(X^*)=\infty$.
Since
$$
\sum_{n=1}^\infty\frac{\nu_n^*(B_n\triangle S_g^*B_n)}{\nu_n^*(B_n)}=
2e\sum_{n=1}^\infty(\nu_n^*(B_n)-\nu_n^*(B_n\cap S_g^*B_n))=
2\sum_{n=1}^\infty(1-e^{1-\nu_n(F_n\cap S_gF_n)})
$$
and
$$
\sum_{n=1}^\infty (1-\nu_n(F_n\cap S_gF_n))=\frac12\sum_{n=1}^\infty
\frac{\kappa(F_n\triangle S_gF_n)}{\kappa(F_n)},
$$
 it follows  from this and $(\alpha_3)$ that 
$$
\sup_{g\in K}\sum_{n=1}^\infty\frac{\nu_n^*(B_n\triangle S_g^*B_n)}{\nu_n^*(B_n)}<\infty\qquad\text{ for each compact $K\subset G$.}\tag3-1
$$
Therefore, by Fact~1.11, $T^*$ preserves $(\nu^*)^{\boldkey B}$.

Thus,  $(X^*,\goth B^*,\mu^*,T^*)$ admits an equivalent invariant infinite $\sigma$-finite measure $(\nu^*)^{\boldkey B}$.
Since $S$ is free, $S$ is effective.
Therefore Proposition~2.4 and~$(\alpha_5)$ yield that $(X^*,\goth B^*,\mu^*,T^*)$ is free.

Unfortunately, we do not know whether  $(X^*,\goth B^*,\mu^*,T^*)$ is ergodic or not.
However we observe that for each increasing sequence of positive integers $(k_n)_{n=1}^\infty$, the conditions  $(\alpha_1)$--$(\alpha_5)$ hold if we replace $\boldkey F$, $(a_n)_{n=1}^\infty$ and $(\lambda_n)_{n=1}^\infty$ with the subsequences
$\boldkey F':=(F_{k_n})_{n=1}^\infty$, $(a_{k_n})_{n=1}^\infty$
and $(\lambda_{k_n})_{n=1}^\infty$ respectively.
Hence given $(Y,\goth Y,\kappa, S)$ and these 3 subsequences, the general construction from \S 2 yields a new dynamical system which we denote by $(\widetilde X,\widetilde{\goth B},\widetilde\mu,\widetilde T)$.
As we have shown above, $(\alpha_1)$--$(\alpha_5)$ imply that the nonsingular Poisson suspension $(\widetilde X^*,{\widetilde{\goth B}}^*,\widetilde\mu^*,\widetilde T^*)$ of  $(\widetilde X,\widetilde{\goth B},\widetilde\mu,\widetilde T)$ is well defined. 
Moreover, the dynamical system $(\widetilde X^*,{\widetilde{\goth B}}^*,\widetilde\mu^*,\widetilde T^*)$
is free and  IDPFT.
Of course, $(\widetilde X^*,{\widetilde{\goth B}}^*,\widetilde\mu^*,\widetilde T^*)$
is canonically isomorphic to $\bigotimes_{n\in\Bbb N}(Y^*,\goth Y^*,\mu_{k_n}^*,S^*)$.
If we let $\boldsymbol B':=(B_{k_n})_{n=1}^\infty$ then
the measure $(\nu^*)^{\boldsymbol B'}$ is infinite, $\sigma$-finite, $\widetilde\mu^*$-equivalent and $\widetilde T^*$-invariant.

\comment
if one replaces $\boldkey F$ with an arbitrary  sub-sequence $\boldkey F':=(F_{k_n})_{n=1}^\infty$ and considers subsequences $(a_{k_n})_{n=1}^\infty$
and $(\lambda_{k_n})_{n=1}^\infty$ instead of $(a_n)_{n=1}^\infty$ and $(\lambda_n)_{n=1}^\infty$ respectively then  then $\boldkey F'$ is also $S$-F{\o}lner and
 $(\alpha_2)$--$(\alpha_5)$ hold for the triplet $(F_{k_n},a_{k_n},\lambda_{k_n})_{n=1}^\infty$.
 Hence, given

 and associates a $G$-action $T'$ with $B$ in the same way as we associated $T$ with $B$ then
$T'$ possesses the same properties that we established for $T$: it is free,  IDPFT and it admits an exhaustingT?-F¿lnersequence. Therefore, to complete the proof of (i)?(ii), it suffices 
\endcomment

{\sl Claim B.} There is an increasing sequence $(k_n)_{n=1}^\infty$ of  positive integers such that  $\bigotimes_{n\in\Bbb N}(Y^*,\goth Y^*,\mu_{k_n}^*,S^*)$  is of infinite ergodic index.

To simplify the argument, we only show how to choose  $(k_n)_{n=1}^\infty$ so that the Cartesian square of $\bigotimes_{n\in\Bbb N}(Y^*,\goth Y^*,\mu_{k_n}^*,S^*)$ is ergodic. 
The sequence  $(k_n)_{n=1}^\infty$  will be defined inductively.
First, we set $k_1:=1$.
Suppose that we have already specified $(k_j)_{j=1}^{n}$.
Our purpose is to define $k_{n+1}$.
For an arbitrary $j\in\Bbb N$, let $U_{j,S^*}$ denote the Koopman unitary representation of $G$ associated with the probability preserving system $(Y^*,\nu_j^*, S^*)$.
We deduce from~$(\alpha_2)$ and Facts~1.2 and~1.21 that the  unitary representation
$U_{j,S^*}\restriction L^2_0(Y^*,\nu_j^*)$ of $G$ is weakly mixing.
Hence
 the dynamical system  $(Y^*,\nu_j^*, S^*)$ is weakly mixing by Fact~1.4(ii).
 It follows that 
 \roster
 \item"$(\triangleright)$" the dynamical system $\bigotimes_{l=1}^{n}(Y^*\times Y^*,\nu_{k_l}^*\otimes \nu_{k_l}^*, (S^*_g\times S^*_g)_{g\in G})$ is ergodic.
 \endroster
   Let $(\Cal P_{j,l})_{l=1}^\infty$ be a  sequence of finite  partitions of $Y^*$ into Borel subsets of equal measure $\nu_j^*$
  such that 
  \roster
   \item"$(\circ)$" $(\Cal P_{j,l})_{l=1}^\infty$ approximates $(\goth Y^*,\nu_j^*)$.
  \endroster
  It follows from 
$(\triangleright)$ that  there is a finite subset $H_n\subset G$ such that for every two atoms $P,Q$
of the partition $\bigotimes_{j=1}^{n}(\Cal P_{k_j,n})^{\otimes 2}$ of $(Y^*\times Y^*)^n$, there is a family of measured subsets $(P_f)_{f\in H_n}$ of $P$ such that
\roster
 \item"$(i)$"
  $P_f\cap P_h=\emptyset$ and $(S_f^*)^{\times 2n}P_f\cap (S_f^*)^{\times 2n}P_h=\emptyset$ for all $f,h\in H$ with  $f\ne h$,
  \item"$(ii)$" $\bigsqcup_{f\in F}(S_f^*)^{\times 2n}P_f\subset Q$ and
  \item"$(iii)$" $\Big(\bigotimes_{j=1}^n(\nu_{k_j}^*)^{\otimes 2}\Big)(\bigsqcup_{f\in H_n}P_f)>0.5\Big(\bigotimes_{j=1}^n(\nu_{k_j}^*)^{\otimes 2}\Big)(P)$.
\endroster
Utilizing $\thetag{3-1}$ we can select $k_{n+1}\in\Bbb N$ large so  that $k_{n+1}>k_n$ and
 $$
  \min_{f\in H_n}\prod_{j\ge k_{n+1}}\frac{\nu_j^*(S_f^*B_j\cap B_j)}{\nu_j^*(B_j)}>0.5.\tag3-2
  $$
Repeating this process infinitely many times, we construct the entire sequence $(k_n)_{n=1}^\infty$.
\comment

Obviously, $(\alpha_1)$--$(\alpha_4)$ hold for the subsequences $(F_{k_n})_{n=1}^\infty$, $(a_{k_n})_{n=1}^\infty$ and $(\lambda_{k_n})_{n=1}^\infty$ of  $(F_{n})_{n=1}^\infty$, $(a_{n})_{n=1}^\infty$ and $(\lambda_{n})_{n=1}^\infty$ respectively.
Hence given $(Y,\goth Y,\kappa, S)$ and these 3 subsequences, the general construction from \S 2 yields a new dynamical system which we denote by $(\widetilde X,\widetilde{\goth B},\widetilde\mu,\widetilde T)$.
The nonsingular Poisson suspension $(\widetilde X^*,{\widetilde{\goth B}}^*,\widetilde\mu^*,\widetilde T^*)$ of  $(\widetilde X,\widetilde{\goth B},\widetilde\mu,\widetilde T)$ is well defined. 
We note that  the dynamical system $(\widetilde X^*,\widetilde{\goth B}^*,\widetilde\mu^*,\widetilde T^*)$  is nothing but the infinite direct product $\bigotimes_{n\in\Bbb N}(Y^*,\goth Y^*,\mu_{k_n}^*,S^*)$.
Thus, $(\widetilde X^*,\widetilde{\goth B}^*,\widetilde\mu^*,\widetilde T^*)$  is
 IDPFT.
If we let  $\widetilde{\boldkey B}:=(B_{k_n})_{n=1}^\infty$
then $\widetilde T^*$ admits an equivalent infinite $\sigma$-finite invariant measure $(\widetilde\nu^*)^{\widetilde{\boldkey B}}$, where $\widetilde\nu^*:=(\nu_{k_n}^*)_{n=1}^\infty$.
We now verify that the dynamical system $(\widetilde X,\widetilde{\goth B},(\widetilde\nu^*)^{\widetilde{\boldkey B}},
\widetilde T)$ is weakly mixing.

\endcomment
Since $\widetilde\mu^*\sim (\nu^*)^{\boldkey B'}$, it suffices to 
verify that the Cartesian square of  $(\widetilde X^*,\widetilde{\goth B}^*,(\nu^*)^{{\boldkey B}'},
\widetilde T^*)$ is ergodic.
As $\widetilde T^*$ preserves $(\nu^*)^{{\boldkey B}'}$, we will use Lemma~3.3 for 
this verification.
Let
   $$
   \widetilde {\boldkey B}^n:=(Y^*)^n\times B_{k_{n+1}}\times B_{k_{n+2}}\times\cdots\subset \widetilde X^*\qquad\text{   for each $n\in\Bbb N$.}\tag3-3
   $$ 
  For each pair of atoms
$P,Q\in\bigotimes_{j=1}^{n}(\Cal P_{k_j,n})^{\otimes 2}$ and $f\in H_n$,
we  set
$$
\align
 P' &:=P\times (B_{k_{n+1}}\times B_{k_{n+2}}\times\cdots)^{\times 2}\in (\Cal P_n')^{\otimes 2},\\
 Q' &:=Q\times(B_{k_{n+1}}\times B_{k_{n+2}}\times\cdots)^{\times 2}\in (\Cal P_n')^{\otimes 2}
 \qquad \text{and}\\ 
 P_f' &:=P_f\times\big((B_{k_{n+1}}\cap (S^*_f)^{-1}B_{k_{n+1}})\times (B_{k_{n+2}}\cap (S^*_f)^{-1}B_{k_{n+2}})\times\cdots\big)^{\times 2},
 \endalign
 $$
 where $(P_f)_{f\in H_n}$ are mutually disjoint subsets of $P$ satisfying $(i)$--$(iii)$.
 Then $\Cal P_n':=(P')_{P\in\bigotimes_{j=1}^{n}(\Cal P_{k_j,n})^{\otimes 2}}$ is a finite partition of $\widetilde{\boldkey B}^n\times\widetilde{\boldkey B}^n$ into subsets of equal measure $(\nu^*)^{{\boldkey B}'}\otimes (\nu^*)^{{\boldkey B}'}$.
 Moreover, $(\circ)$  yield that
 \roster
 \item"$(\bullet)$"
%
$(\Cal P_n')_{n=1}^\infty$ approximates
$((\widetilde{\goth B}^*\otimes\widetilde{\goth B}^*)_0,(\nu^*)^{{\boldkey B}'}\otimes (\nu^*)^{{\boldkey B}'})$,
   \endroster
   where $(\widetilde{\goth B}^*\otimes\widetilde{\goth B}^*)_0$ denote the ring of subsets
   of finite measure in $\widetilde{\goth B}^*\otimes\widetilde{\goth B}^*$.
 We also note that
$(P'_f)_{f\in H_n}$ are  mutually disjoint subsets of $P'$
and $(( \widetilde T_f^*\times \widetilde T_f^*)P'_f)_{f\in H_n}$ 
are  mutually disjoint subsets of $Q'$.
It follows from \thetag{3-2} that for each $f\in H_n$,
$$
\prod_{j\ge n+1}\frac{\nu_{k_j}^*(S_f^*B_{k_j}\cap B_{k_j})}{\nu_{k_j}^*(B_{k_j})}>0.5.
  $$
Using this inequality and $(iii)$, we obtain that
$$
\align
\Big((\nu^*)^{{\boldkey B}'}\otimes (\nu^*)^{{\boldkey B}'}\Big)\bigg(\bigsqcup_{f\in H_n}P'_f\bigg)&>
\frac{\frac12\Big(\bigotimes_{j=1}^n(\nu_{k_j}^*)^{\otimes 2}\Big)(P)}{\Big(\prod_{j=1}^n\nu^*_{k_j}(B_{k_j})\Big)^2}\cdot \frac14
\\
&=\frac 18
\Big((\nu^*)^{{\boldkey B}'}\otimes (\nu^*)^{{\boldkey B}'}\Big)(P')
\endalign
$$
for every pair of $\Cal P_n'$-atoms $P',Q'$.
Therefore it follows from Lemma~3.3 that the $G$-action  $(\widetilde T_g^*\times \widetilde T_g^*)_{g\in G}$ on 
$(\widetilde{X}^*\otimes\widetilde{X}^*,(\nu^*)^{{\boldkey B}'}\otimes (\nu^*)^{{\boldkey B}'})$ is ergodic.
Thus, Claim~B is proved.

It follows that $(\widetilde X^*,\widetilde{\goth B}^*,\widetilde\mu^*,\widetilde T^*)$ is of type $II_\infty$.
We deduce from Proposition~2.3 and $(\alpha_1)$ that $(\widetilde X^*,\widetilde{\goth B}^*,\widetilde\mu^*,\widetilde T^*)$
is not strongly ergodic.
Since $(\widetilde X^*,\widetilde{\goth B}^*,\widetilde\mu^*,\widetilde T^*)$ is IDPFT,  Proposition~1.25 yields that $T^*$ is amenable in the Greenleaf sense.
Moreover, $\widetilde T^*$ is amenable if and only if $G$ is amenable.

Thus, Theorem~3.2 is proved completely. \qed
\enddemo


\comment

\remark{Remark \rom{3.3}}  It is possible to strengthen the claim on  weak mixing of 
$(\widetilde X^*,\widetilde \mu^*,\widetilde T^*)$
in  Theorem~3.1.
Passing to a further subsequence in $(k_n)_{n=1}^\infty$, we can construct  $(\widetilde X,\widetilde{\goth B},\widetilde\mu, \widetilde T)$ in such a way that the system $(\widetilde X^*,\widetilde{\goth B}^*, \widetilde\mu^*, \widetilde T^*)$ is of infinite ergodic index.
We leave details to the reader.
\endremark

\endcomment

\comment

\remark{Remark \rom{3.4}} There is a $\widetilde T^*$-F{\o}lner sequence $(C_n)_{n=1}^\infty$
such that $(\nu^*)^{{\boldkey B}'}(C_n)=1$ for all $n\in\Bbb N$.
This assertion
  strengthens the claim on the Greenleaf's amenability of $(\widetilde X^*, \widetilde\mu^*, \widetilde T^*)$ in the statement of Theorem~3.2.
First, we note that the sequence $(\widetilde {\boldkey B}^n)_{n=1}^\infty$ defined in \thetag{3-3} is a $\widetilde T^*$-F{\o}lner sequence
 (see a remark just below Fact~1.11).
 However   $\widetilde {\boldkey B}^1\subset\widetilde {\boldkey B}^2\subset\cdots$ with $\bigcup_{n=1}^\infty\widetilde {\boldkey B}^n=\widetilde X^*$.
Now, suppose
that
\roster
\item"$(\alpha_6)$" $\sup_{g\in K}\frac{n\kappa(F_{k_n}\triangle S_gF_{k_n})}{\kappa(F_{k_n})}\to 0$ as $n\to \infty$ for each compact subset $K\subset G$ and
\item"$(\alpha_7)$" there is 
 a sequence $(A_n)_{n=1}^\infty$ of subsets  from $\goth Y$ such that $\kappa(A_n)=n\kappa (F_{k_n})$
 and $\kappa(A_n\triangle S_gA_n)=n\kappa(F_{k_n}\triangle S_gF_{k_n})$
  for each $n\in\Bbb N$.
 \endroster
 Then we define a sequence $(C_n)_{n=1}^\infty$ of subsets in $X^*$ by setting
$$
C_n:=(Y^*)^{n-1}\times [A_{n}]_0\times B_{k_n+1}\times B_{k_n+2}\times\cdots.
$$
We now have that $\nu_{k_n}(A_n)=n$
and
$$
(\widetilde\nu^*)^{\widetilde{\boldkey B}}(C_n)=\frac{\nu^*_{k_n}([A_{n}]_0)}{\prod_{j=1}^n\nu^*_{k_j}(B_{k_j})}=
e^{-\nu_{k_n}(A_{n})+\sum_{j=1}^n\nu_{k_j}(F_{k_j})}=1.
$$
On the other hand,
$$
\align
(\widetilde\nu^*)^{\widetilde{\boldkey B}}(C_n\cap \widetilde T^*_gC_n)&=
\frac{\nu^*_{k_n}([A_{n}\cup S_g A_{n}]_0)}{\prod_{j=1}^n \nu^*_{k_j}(B_{k_j})}
\prod_{j>n}\frac{\nu^*_{k_j}(B_{k_j}\cap S_g^*B_{k_j})}{\nu^*_{k_j}(B_{k_j})}\\
&=
\frac{e^{-\nu_{k_n}(A_{n}\cup S_gA_{n})}}{e^{-n}}
e^{-\sum_{j>n}(\nu_{k_j}(F_{k_j}\cup S_gF_{k_j})-1)}\\
&=e^{-0.5 \nu_{k_n}(A_{n}\triangle S_gA_{n})}e^{-0.5\sum_{j>n}
\nu_{k_j}(F_{k_j}\triangle S_gF_{k_j})
}.
\endalign
$$
Hence 
$$
\log \Big((\widetilde\nu^*)^{\widetilde{\boldkey B}}(C_n\cap\widetilde T^*_gC_n)\Big)=
-\frac{\kappa(A_{n}\triangle S_gA_{n})}{2\kappa(F_{k_n})}-\sum_{j>n}
\frac{\kappa(F_{k_j}\triangle S_gF_{k_j})}{2\kappa(F_{k_j})}\to 0
$$
as $n\to\infty$
 in view of  $(\alpha_3)$, $(\alpha_5)$
 and $(\alpha_6)$. 
 Thus,  $(C_n)_{n=1}^\infty$ is a $\widetilde T^*$-F{\o}lner sequence in $\big({\widetilde X}^*,
 \widetilde{ \goth B}^*, 
 (\widetilde\nu^*)^{\widetilde{\boldkey B}}\big)$.

 It remains to show $(\alpha_5)$ and $(\alpha_6)$.
 As for $(\alpha_5)$, it holds on  a subsequence of $(F_{k_n})_{n=1}^\infty$.
  Hence it is enough to drop to this subsequence.
 Next,  
 the simplest way to construct $(A_n)_{n=1}^\infty$  satisfying $(\alpha_6)$, 
is to replace the original dynamical system $(Y,\kappa, S)$ with the infinite disjoint union  $\bigsqcup_{n=1}^\infty (Y,\kappa, S)$.
Then the role of $(F_{k_n})_{n=1}^\infty$ in the new system will be played by the sequence $(F_{k_n})_{n=1}^\infty$ contained in the first copy of $Y$ in $X$.
We now define $A_n$ as the union of $n$ copies of $F_{k_n}$ contained in $n$ different copies of $Y$ in $X$, $n\in\Bbb N$.
Then,  $(\alpha_1)$--$(\alpha_6)$ are all satisfied for 
$\bigsqcup_{n=1}^\infty (Y,\kappa, S)$ and the new $(F_{k_n})_{n=1}^\infty$.
\endremark

\endcomment

\proclaim{Theorem 3.5}
If $G$ is Haagerup then we can choose the  system $(\widetilde X,\widetilde{\goth B},\widetilde\mu, \widetilde T)$ in the statement of  Theorem~3.2 in such a way that the Poisson suspension $(\widetilde X^*,\widetilde{\goth B}^*, \widetilde\mu^*, \widetilde T^*)$   is 
of  0-type (in addition to the other properties of
$(\widetilde X^*,\widetilde{\goth B}^*, \widetilde\mu^*, \widetilde T^*)$ listed in  Theorem~3.2).
\endproclaim

\demo{Proof}
In view of  \cite{DeJoZu} (see also \cite{Da}), there is an infinite $\sigma$-finite measure preserving  system $(Y,\goth Y,\kappa, S)$ of 0-type satisfying $(\alpha_1)$--$(\alpha_5)$.
Then we repeat the  proof of Theorem~3.2 verbally to  obtain a
non-strongly ergodic free  IDPFT system $(\widetilde X^*,\widetilde{\goth B}^*,\widetilde\mu^*,\widetilde T^*)$ of infinite ergodic index and of type $II_\infty$.
The  $\widetilde\mu^*$-equivalent infinite $\sigma$-finite $\widetilde T^*$-invariant measure is $(\nu^*)^{{\boldkey B}'}$.
Since the dynamical system $(Y,\nu_n,S)$ is of 0-type, the Poisson suspension 
 $(Y^*,\nu_n^*,S^*)$ is mixing in view of Fact~1.21 for each $n\in\Bbb N$.
 Then by Corollary~1.13, $(\widetilde X^*,\widetilde{\goth B}^*,(\nu^*)^{{\boldkey B}'},\widetilde T^*)$
is of 0-type.
Hence $(\widetilde X^*,\widetilde{\goth B}^*,\widetilde\mu^*,\widetilde T^*)$
is of 0-type too.
\qed
\enddemo

\head 4. Type $III_0$ ergodic Poisson suspensions
\endhead

Our purpose in this section is to prove 
 the implications (1)$\Rightarrow$(3) from Theorems~A and B for  $K=III_0$.
Thus, we assume that  $G$ does not have property (T).
As in \S3, we 
fix a measure preserving free $G$-action
$S=(S_g)_{g\in G}$ on     an infinite $\sigma$-finite standard measure space $(Y,\goth Y,\kappa)$ such that
\roster
\item"$(\alpha_1)$"
there is an $S$-F{\o}lner sequence $\boldkey F=(F_n)_{n=1}^\infty$ such that $\kappa(F_n)=1$ for each $n\in\Bbb N$,
\item"$(\alpha_2)$" the corresponding unitary Koopman representation $U_S$ of $G$ is weakly mixing,
\item"$(\alpha_3)$" $\sum_{n=1}^\infty\sup_{g\in K}{\kappa(S_gF_n\triangle F_n)}<+\infty$ for each compact $K\subset G$.
\item"$(\alpha_4)$" There is a $\kappa$-preserving transformation $Q$  that commutes with $S_g$ for each $g\in G$ and such that each $Q$-invariant subset is either of 0-measure or of infinite measure.
\endroster
We deduce  from $(\alpha_2)$ and Fact~1.2 that 
there is a  dispersed  sequence $(g_n)_{n=1}^\infty$ in $G$ such that
\roster
\item"$(\alpha_5)$"
$S$ is of $0$-type  along the subset $H:=\{g_ng_m^{-1}\mid n>m, n,m\in\Bbb N\}\subset G$.
\endroster
 Utilizing $(\alpha_3)$ and passing to a suitable subsequence in $\boldkey F$, which we denote by the same symbol $\boldkey F$, we may (and will) assume without loss of generality 
  and 
$$
\sum_{n=1}^\infty\max_{1\le l\le n}{\kappa(S_{g_l}F_n\triangle F_n)} \le 2.\tag{4-1}
$$

\comment

Fix two sequences $(a_n)_{n=1}^\infty$ and $(\lambda_n)_{n=1}^\infty$
of positive reals such that
\roster
\item"$(\alpha_4)$" $a_n:=1$ for all $n\in\Bbb N$ and  $\sum_{n=1}^\infty\lambda_n<+\infty$.
\endroster

In  \cite{BeHaVa, Theorem~6.3.4} (see also \cite{CoWe}), given an arbitrary weakly mixing orthogonal representation $\pi$ of $G$ in a real Hilbert space $\Cal H$ with an almost invariant vector ($\pi$ does not have an invariant vector because $\pi$ is weakly mixing), 
a weakly mixing non-strongly ergodic probability preserving Gaussian $G$-action is constructed.
We introduce another real Hilbert space $\widetilde{\Cal H}:=\bigoplus_{n\in\Bbb Z}\Cal H$
and consider an orthogonal  representation $\widetilde\pi:=\bigoplus_{n\in\Bbb Z}\pi$ of $G$ in $\widetilde{\Cal H}$.
Let $\vartheta$ stand for the left shift on $\widetilde{\Cal H}$, i.e. $\vartheta((h_n)_{n=1}^\infty):=(h_{n+1})_{n\in\Bbb Z}$ for each vector $(h_{n+1})_{n\in\Bbb Z}\in\widetilde{\Cal H}$.
Then $\vartheta$ is  an orthogonal operator that commutes with $\widetilde\pi$.
Of course, $\widetilde\pi$ is weakly mixing and $\widetilde\pi$ has  almost invariant vectors.
Let $R=(R_g)_{g\in G}$ be the measure preserving Gaussian action on the corresponding standard probability space $(Z,\goth Z,\tau)$ constructed as in \cite{BeHaVa, Theorem~6.3.4}
but associated with $\widetilde\pi$.
Then $R$ is non-strongly ergodic weakly mixing.

 Hence
 there is an $R$-mixing sequence $(g_n)_{n=1}^\infty$ in $G$.
 Passing, if necessary, to a subsequence, we may (and will) assume without loss of generality that 
\roster
\item"$(\circ)$"
the sequence $(g_n)_{n=1}^\infty$ is dispersed and
\item"$(\diamond)$"
$R$ is mixing along the subset $H:=\{g_ng_m^{-1}\mid n>m, n,m\in\Bbb N\}\subset G$.
\endroster

Since $R$ is not strongly ergodic,
 there is an $R$-asymptotically invariant sequence $\boldkey P=(P_n)_{n=1}^\infty$
of subsets in $Z$ with $\tau(P_n)=0.5$ for each $n\in\Bbb N$ (see the proof of  \cite{BeHaVa, Theorem~6.3.4} and \cite{CoWe}).
Passing to a subsequence in $\boldkey P$, if necessary, we will assume also that 
\roster
\item"$(\alpha_1)$"
$\sum_{n=1}^\infty\sup_{g\in K}\tau(P_n\triangle R_gP_n)<\infty$
for each compact $K\subset  G$ 
\endroster
Also, a measure preserving transformation $W$ of $(Z,\goth Z,\tau)$ is generated by $\vartheta$. 
 Since the maximal spectral type of $\vartheta$ is continuous, $W$ is weakly mixing.
 Of course, $WR_g=R_gW$  for all  $g\in G$.

Let $(Y,\goth Y):=(Z,\goth Z)^{\otimes \Bbb N}$, $S_g:=R_g^{\times\Bbb N}$ for each $g\in G$ and let $\kappa$ denote the  restricted infinite product of the constant sequence of probabilities $(\tau)_{n=1}^\infty$ with respect to $\boldkey P$.
Since $\prod_{n=1}^\infty\tau(P_n)=0$, it follows that $\kappa$ is infinite and $\sigma$-finite. Moreover,  $(\alpha_1)$ and  Fact~1.11 yield that $S$ preserves $\kappa$.
It follows from Corollary~1.13 and $(\diamond)$ that $S$ is of 0-type  along $H$.
For each $n\in\Bbb N$, we let $ F_n:=Z^n\times P_{n+1}\times P_{n+2}\times\cdots$.
From $(\alpha_1)$ we deduce that 
  the sequence $\boldkey F:=(F_n)_{n=1}^\infty$ is $S$-F{\o}lner (see \thetag{1-2} and the sentence below it).

For each $n>0$, we define a transformation $V_n$ of $Y$ by setting
$$
V_n\big((z_k)_{k=1}^\infty\big):=(z_k')_{k=1}^\infty,
$$
where $z_n':=Wz_n$ and $z_k':=z_k$ if $k\ne n$.
It is straightforward to verify that  $V_n$ preserves $\kappa$.
Denote by $\Gamma$ the  transformation group of $(Y,\goth Y,\kappa)$ generated by all $V_n$,  $n\in\Bbb N$.
Since $W$ commutes with $R$, we obtain that $\Gamma$ commutes with $S$.
Since $W$ is ergodic, it follows that $\Gamma$ is ergodic.

\endcomment

Let $(l_n)_{n=1}^\infty$ be  a sequence of positive integers such that $l_1|l_2$, $l_2|l_3,\dots$, $\lim_{n\to\infty} l_n=+\infty$ and
$\sum_{n=1}^\infty\frac1{n4^{l_n}}=+\infty$.
We now set
\roster
\item"$(\alpha_6)$" $\lambda_n:=2^{l_n}$ and $a_n:=c(\lambda_n)^{-1}n^{-1}$ for each $n>0$.
\endroster
Denote by $(X,\goth B,\mu, T)$ the dynamical system associated with $(Y,\goth Y,\kappa,S)$,  $\boldkey F$, $(a_n)_{n=1}^\infty$ and $(\lambda_n)_{n=1}^\infty$ via the general construction from \S 2.
Since $a_n\sim\lambda_n^{-3}n^{-1}$ as $n\to\infty$, it follows from $(\alpha_3)$ and $(\alpha_6)$ that for each $g\in G$,
$$
\gather
\sum_{n\in\Bbb N}a_n|\lambda_n-1|\frac{\kappa
(S_gF_n\triangle F_n)}{\kappa(F_n)}<\infty\qquad\text{and}\tag4-2\\
\sum_{n=1}^\infty 
a_n\frac{
|1-\lambda_n^2|(1+\lambda_n^3)}{\lambda_n^2}
\frac{\kappa(S_gF_n\triangle F_n)}{\kappa(F_n)}<\infty.\tag4-3
\endgather
$$
Hence \thetag{4-2} and Lemma~2.1 yield that $T_g\in\text{Aut}_1(X,\mu)$  and $\chi(T_g)=0$ for each $g\in G$.
Therefore, as was explained in \S2,  the nonsingular Poisson suspension $(X^*,\goth B^*,\mu^*,T^*)$ of 
$(X,\goth B,\mu, T)$ is well defined.

\proclaim{Theorem~4.1} The dynamical system $(X^*,\goth B^*,\mu^*,T^*)$ is weakly mixing non-strongly ergodic  IDPFT of Krieger type $III_0$.
Hence $T^*$ is amenable in the Greenleaf sense.
Moreover,  $T^*$ is amenable if and only if $G$ is amenable.
\endproclaim
\demo{Proof}
It was shown in \S2 that $(X^*,\goth B^*,\mu^*,T^*)$ is IDPFT.
Hence $T^*$ is amenable in the Greenleaf sense by Proposition~1.25.
Proposition~2.4 and $(\alpha_4)$ imply that $T^*$ is free.
 Proposition~2.3 and $(\alpha_1)$ imply that  $T^*$ is not strongly ergodic.

We now prove that $T^*$ is weakly mixing.
We note that  \thetag{4-3}, Lemma~2.2  and~$(\alpha_6)$ yield that
 $$
 \gather
 \bigg(\frac {d\mu}{d\mu\circ T_g^{-1}}\bigg)^2-1\in L^1(X,\mu)\qquad\text{and}\\
\int_X\Big(\Big(\frac{d\mu}{d\mu\circ T_g^{-1}}\Big)^2-1\Big)d\mu=\sum_{n=1}^\infty a_nc(\lambda_n)\frac{\kappa
(S_gF_n\triangle F_n) }{2\kappa(F_n)}=\sum_{n=1}^\infty\frac1{2n}{\kappa
(S_gF_n\triangle F_n) }
\endgather
$$
for every $g\in G$.
Hence, for each $k\in\Bbb N$,
$$
\align
\int_X\Big(\Big(\frac{d\mu}{d\mu\circ T_{g_k}^{-1}}\Big)^2-1\Big)d\mu
&\le\frac 12\sum_{n=1}^k\frac 1n+ 
\sum_{n=k+1}^\infty{\kappa
( S_{g_k} F_n\triangle F_n )}\le \frac 12\log k+3
\endalign
$$
in view of \thetag{4-1}.
This implies that
$$
\sum_{k=1}^\infty \frac 1{k^2}\,e^{\int_X\Big(\Big(\frac{d\mu}{d\mu\circ T_{g_k}^{-1}}\Big)^2-1\Big)d\mu}<+\infty.
$$
Therefore, by Proposition 1.20,  
 $$
 \sum_{k=1}^\infty\frac{d\mu^*\circ (T_{g_k}^*)^{-1}}{d\mu^*}(\omega)=+\infty\qquad\text{ at a.e. $\omega$.}
$$
On the other hand, in view of  Fact~1.21,
 $(Y^*,\goth Y^*,\nu_k^*, S^*)$ is mixing  along $H$  for each $k\in \Bbb N$ because $(\alpha_5)$ holds.
 (We refer to \S2 for the  definitions of $(\nu_k)_{k=1}^\infty$ and $(\mu_k)_{k=1}^\infty$ .)
Hence by Proposition~1.16, $T^*$ is weakly mixing.

{\sl Claim A.} The dynamical system $(X^*,\goth B^*,\mu^*,T^*)$ is of type $III$.

Suppose, by contraposition, that $T^*$ is of type $II$.
Let $\eta$ stand for a $T^*$-invariant $\sigma$-finite measure equivalent to $\mu^*$.
Fix $n\in\Bbb N$.
Since $\bigotimes_{k=1}^n\mu_k^*\sim\bigotimes_{k=1}^n\nu_k^*$, it follows that the projection of $\eta$ to $(Y^*)^n$ 
along the mapping 
$$
\pi_n:X^*\ni(x^*_k)_{k=1}^\infty\mapsto (x^*_k)_{k=1}^n\in( X^*)^n
$$
 has the same collection of 0-measure subsets as the measure $\bigotimes_{k=1}^n\nu_k^*$ has.
 Hence we can disintegrate $\eta$ over $\bigotimes_{k=1}^n\nu_k^*$ along $\pi_n$.
 Thus, there is measurable field $(Y^*)^n\ni y\mapsto \eta_y$ of $\sigma$-finite measures on $(Y^*)^{\Bbb N}$ such that
$$
\eta=\int_{(Y^*)^n}\delta_y\otimes\eta_y\,d\bigg(\bigotimes_{k=1}^n\nu_k^*\bigg)(y).\tag4-4
$$
Since $Q$ preserves $\nu_k$, we obtain that the Poisson suspension $Q^*$
of $Q$ preserves $\nu_k^*$.
It follows from $(\alpha_4)$ that $Q^*$ is ergodic
with respect to $\nu_k^*$ for each $k\in\Bbb N$.
Moreover, $Q^*$ commutes with $S^*$.
Let 
$$
\Gamma_n:=\{(Q^*)^{l_1}\otimes\cdots\otimes  (Q^*)^{l_n}\mid l_1,\dots,l_n\in\Bbb Z\}.
$$
Then $\Gamma_n$ is a countable group (isomorphic to $\Bbb Z^n$) of measure preserving transformations of the product space
$((Y^*)^n,(\goth Y^*)^{\otimes n},\bigotimes_{k=1}^n\nu_k^*)$.
Of course,  $\Gamma_n\otimes 1$ is a group of $\mu^*$-nonsingular transformations commuting with $T^*$.
Since $\eta\sim\mu^*$,  $\Gamma_n\otimes 1$ is also $\eta$-nonsingular.
Since $T^*$ is ergodic and $\eta$-preserving, it follows that 
$\eta\circ \gamma=a_\gamma\eta$ for some $a_{\gamma}>0$, for each transformation $\gamma\in\Gamma_n\otimes 1$.
Since $\gamma$ preserves a finite measure $\Big(\bigotimes_{k=1}^n\nu_k^*\Big)\otimes\bigotimes_{k>n}\mu_n^*$ which is equivalent to $\eta$, we obtain that $\gamma$ is conservative with respect to $\eta$.
Hence $a_\gamma=1$, i.e. $\gamma$ preserves $\eta$ for each $\gamma\in \Gamma_n\otimes 1$.
We deduce from this fact and \thetag{4-4} that the mapping 
$$
\Psi:(Y^*)^n\ni y\mapsto \eta_y
$$
 is invariant under $\Gamma_n$.
 Since $\Gamma_n$  is ergodic with respect to  $\bigotimes_{k=1}^n\nu_k^*$, it follows that
$\Psi$ is constant almost everywhere, i.e.
 there  exists a $\sigma$-finite measure $\eta^n$ on $((Y^*)^\Bbb N,(\goth Y^*)^{\otimes \Bbb N}),$ such that
$\eta_y=\eta^n$ at a.e. $y\in (Y^*)^n$.
We thus obtain that 
 $\eta=\big(\bigotimes_{k=1}^n\nu_k^*\big)\otimes\eta^n$.
 Since   $n$ is arbitrary,
 $\eta$ is an MH-product of $(\nu_k^*)_{k=1}^\infty$ (see Definition~1.10(ii)).
 Hence, by Fact~1.14(ii), $\eta$ is proportional to the restricted product  of $(\nu^*_k)_{k=1}^\infty$ with respect to some sequence $\boldkey C=(C_n)_{n=1}^\infty$ of subsets  $C_n\subset Y^*$.
 In particular,
 $$
  \eta(C_1\times C_2\times\cdots)>0.\tag4-5
 $$
If $\prod_{n=1}^\infty\nu_n^*(C_n)>0$ then $\eta\sim\bigotimes_{n=1}^\infty\nu_n^*$ and hence
 $\mu^*\sim\nu^*$.
The latter happens if and only if  
$
\sqrt{\frac{d\mu}{d\nu}}-1\in L^2(X,\nu)
$
by Fact~1.18(ii).
However,
$$
\int_X\bigg(\sqrt{\frac{d\mu}{d\nu}}-1\bigg)^2d\nu=\sum_{n=1}^\infty(\sqrt{\lambda_n}-1)^2a_n\asymp \sum_{n=1}^\infty\frac1{n\lambda_n^2}=\sum_{n=1}^\infty\frac 1{n4^{l_n}}=+\infty
$$
by  choice of $(\lambda_n)_{n=1}^\infty$, a contradiction.
Hence $\prod_{n=1}^\infty\nu_n^*(C_n)=0$ or, equivalently,
$\sum_{n=1}^\infty\nu_n^*(Y_n^*\setminus C_n)=+\infty$.
Since  $\eta\sim\mu^*$, we deduce from \thetag{4-5} that
$$
0<\mu^*(C_1\times C_2\times\cdots)=\prod_{n=1}^\infty\mu_n^*(C_n)
$$
 or, equivalently,
 $\sum_{n=1}^\infty\mu_n^*(Y_n^*\setminus C_n)<+\infty$.
 On the other hand, it follows from Fact~1.18(iii)  that
 $$
 \frac{d\nu_n^*}{d\mu_n^*}(\omega)=e^{-\big(\frac1{\lambda_n}-1\big)\mu_n(F_n)}\lambda_n^{-\omega(F_n)}
 $$
at $\mu_n^*$-a.e. $\omega\in Y_n^*$. 
 Since $\lambda_n\ge 1$ and $\mu_n(F_n)=a_n\lambda_n=\frac{\lambda_n}{nc(\lambda_n)}\le 1$, we obtain that 
 $\frac{d\nu_n^*}{d\mu_n^*}(\omega)\le e$ for a.e. $\omega\in X_n^*$.
Hence 
$$
+\infty=\sum_{n=1}^\infty\nu_n^*(Y_n^*\setminus C_n)\le e^{}\sum_{n=1}^\infty\mu_n^*(Y_n^*\setminus C_n)<+\infty,
$$ 
a contradiction.
Thus,  Claim A is proved.

{\it Claim B.} $(X^*,\goth B^*,\mu^*,T^*)$ is of type $III_0$.

Fix $n\in\Bbb N$.
Let $\mu^{*,n}$ denote the infinite product measure $\big(\bigotimes_{k=1}^n\nu_n^*\big)\otimes\bigotimes_{k>n}\mu_k^*$.
Of course,  $\mu^{*,n}\sim\mu^*$.
Hence $r(T^*,\mu^*)=r(T^*,\mu^{*,n})$.
We recall that
$$
f_n:=\frac{d\mu_n}{d\nu_n}=1+(\lambda_n-1)1_{F_n}.
$$
From Fact~1.18(iv) we deduce that
$$
\frac{d\mu_k^*\circ S_g^*}{d\mu_k^*}(\omega)=\prod_{\{y\in Y\mid \omega(\{y\})=1\}}\frac{d\mu_k\circ S_g}{d\mu_k}(y)=
\prod_{\{y\in Y\mid \omega(\{y\})=1\}}\frac{f_k(S_gy)}{f_k(y)}
\in\{\lambda_k^m\mid m\in\Bbb Z\}
$$
for each $g\in G$ at $\mu_k^*$-a.e. $\omega\in Y^*$.
We deduce from this,  $(\alpha_6)$ and the fact that $(X^*,\goth B^*,\mu^{*,n}, T^*)$ is IDPFT that
$$
\frac{d\mu^{*,n}\circ T_g^*}{d\mu^{*,n}}(\omega)=\prod_{k>n}^\infty\frac{d\mu_k^*\circ S_g^*}{d\mu_k^*}(\omega_k)\in\{2^{kl_{n+1}}\mid k\in\Bbb Z\}
$$
for each $g\in G$ at a.e. $\omega\in X^*$.
Hence $r(T^*,\mu^{*,n})\subset\{2^{kl_{n+1}}\mid k\in\Bbb Z\}$.
This yields that $r(T^*,\mu^{*})\subset\bigcap_{n=1}^\infty\{2^{kl_{n+1}}\mid k\in\Bbb Z\}=\{1\}$.
Since $T^*$ is not of type $II$, it follows that $T^*$ is of type $III_0$.
Thus, Claim~B is proved.
\qed

\comment

{\sl Claim C.} 
$(X^*,\goth B^*,\mu^{*}, T^*)$ is not strongly ergodic.

First, we recall that $\widetilde {\Cal H}=\bigoplus_{n\in\Bbb Z}\Cal H$.
Hence $(Z,
\goth Z,\tau, R)=\bigotimes_{n\in\Bbb Z}(Z',\goth Z',\tau', R')$, where $(Z',\goth Z',\tau', R')$
is the Gaussian dynamical system associated with the pair $(\Cal H,\pi)$ as in 
 \cite{BeHaVa, Theorem~6.3.4}.
According to \cite{BeHaVa, Theorem~6.3.4} and \cite{CoWe}, there is  an $R'$-asympto\-tically invariant sequence $(P_n')_{n\in\Bbb Z}$ such that
$\tau'(P_n)=0.5$ for each $n\in\Bbb N$.
For each $m\in\Bbb N$, we let $P_{m,n}:=P_n^{\times m}\times Z'\times Z'\times\cdots\subset Z$.
A straightforward verification shows that the sequence
$(P_{m,n})_{n=1}^\infty$ is $R$-asymptotically invariant and $\tau(P_{m,n})=2^{-m}$ for all $n,m\in\Bbb N$.
Utilizing  the diagonalization argument, we can choose an $R$-asymptotically invariant sequence $(D_n)_{n=1}^\infty$ such that $\tau(D_n)=2^{-n}$ for all $n\in\Bbb N$.
We now let for each $n\in\Bbb N$,
$$
D_n':=Z^{n-1}\times D_n\times P_n\times P_{n+1}\times\cdots\subset  Y.
$$
It is straightforward to check that $\kappa (D_n')=1$ for each $n\in\Bbb N$ and that
$(D_n')_{n=1}^\infty$ is an $S$-F{\o}lner sequence.
Then, by Proposition~2.3, 
$(X^*,\goth B^*,\mu^{*}, T^*)$  is not strongly ergodic.
Thus, Claim~C is proved.

To complete the  proof of Theorem~4.1, it remains to show that 
$T^*$ is amenable in the Greenleaf sense and that 
 $T^*$ is amenable if and only if $G$ is amenable.
This follows from Proposition~1.25.

\endcomment

\enddemo

\remark{Remark \rom{4.2}}  Changing the   parameters $(a_n)_{n=1}^\infty$ and $(l_n)_{n=1}^\infty$ in $(\alpha_6)$ in an appropriate way we can obtain $T^*$ with infinite ergodic index (which is stronger than the  weak mixing).
For that, we choose $(l_n)_{n=1}^\infty$  so that  $\sum_{n=1}^\infty\frac{1}{n4^{l_n}\log(n+1)}=+\infty$.
Then we let
$\lambda_n:=2^{l_n}$ and
$a_n:=c(\lambda_n)^{-1}\frac 1{n\log{(n+1)}}$ for each $n\in\Bbb N$. 
Let $(X,\goth B,\mu,T)$ be the  dynamical system associated with these new parameters
via the general construction from \S2.
Then the Poisson suspension  $(X^*,\goth B^*,\mu^*, T^*)$ possesses all the properties listed in the statement of Theorem~4.1 and, in addition, $T^*$ is of  infinite ergodic index.
Indeed, fix $l>1$.
By Fact~1.18(v),
 the $l$-th Cartesian power  of $(X^*,\goth B^*,\mu^*, T^*)$ is the Poisson suspension of the disjoint union  $\bigsqcup_{j=1}^l(X,\goth B,\mu, T)$ of $l$  copies of $(X,\goth B,\mu, T)$. 
 It is straightforward to verify  that
$\bigsqcup_{j=1}^l(X,\goth B,\mu, T)$
  is the dynamical system associated with 
the disjoint union $\bigsqcup_{j=1}^l(Y,\goth Y,\kappa, S)$ of $l$ copies of $(Y,\goth Y,\kappa, S)$, and the sequences 
$(\bigsqcup_{j=1}^l F_n)_{n=1}^\infty$,  
$(la_n)_{n=1}^\infty$ and $(\lambda_n)_{n=1}^\infty$ via the general construction of \S 2.
Then a
slight modification of the proof of Theorem~4.1 yields that 
$(X^*,\goth B^*,\mu^*, T^*)^{\times l}$ is weakly mixing. 
We leave details to the reader.
Thus, $(X^*,\goth B^*,\mu^*, T^*)$ is of infinite ergodic index.
\endremark


\proclaim{Theorem 4.3}
If $G$ is Haagerup then there is a dynamical system $(X,\goth B,\mu, T)$ so that the Poisson suspension $( X^*,\goth B^*,\mu^*,  T^*)$   is 
of 0-type (in addition to the other properties of $( X^*,\goth B^*,\mu^*,  T^*)$ listed in Theorem~4.1).
\endproclaim

\demo{Proof} Since $G$ is Haagerup, there is a measure preserving $G$-action $S=(S_g)_{g\in G}$ on an infinite $\sigma$-finite standard measure space $(Y,\goth Y,\kappa)$ such that
$(\alpha_1)$--$(\alpha_6)$ hold and, moreover, $S$ is of 0-type (see \cite{DeJoZu} or \cite{Da}).
Then, following the proof of Theorem~4.1 verbally, we construct a dynamical system
$( X,\goth B,\mu,  T)$ whose Poisson suspension $( X^*,\goth B^*,\mu^*,  T^*)$ possesses 
all the properties listed in the statement of Theorem~4.1.
It remains  only to prove that $( X^*,\goth B^*,\mu^*,  T^*)$ is of 0-type.
By \cite{DaKoRo1, Proposition~6.9},
$T^*$ is of 0-type if and only if $\int_X\Big(\sqrt{\frac{d\mu\circ T_g}{d\mu}}-1\Big)^2d\mu\to 0$ as $g\to\infty$.
Since 
$$
\int_X\Bigg(\sqrt{\frac{d\mu\circ T_g}{d\mu}}-1\Bigg)^2d\mu\le
\int_X\Big|{\frac{d\mu\circ T_g}{d\mu}}-1\Big|d\mu
$$
and, as was computed in the proof of Lemma~2.1,
$$
\int_X\Big|{\frac{d\mu\circ T_g}{d\mu}}-1\Big|d\mu=\sum_{n=1}^\infty a_n|\lambda_n-1|\kappa(S_gF_n\triangle F_n),
$$
it suffices to show that $\sum_{n=1}^\infty a_n|\lambda_n-1|(1-\kappa(S_gF_n\cap F_n))\to 0$
as $g\to\infty$.
We recall that $a_n\lambda_n\sim\frac 1{n\lambda_n^2}=\frac1{n4^{l_n}}$ and $\sum_{n=1}^\infty\frac 1{n4^{l_n}}<\infty$.
Hence given $\epsilon>0$, there is $N>0$ such that 
$
 \sum_{n=N+1}^\infty a_n|\lambda_n-1|(1-\kappa(S_gF_n\cap F_n))<\epsilon.
 $
On the other hand, as $S$ is of 0-type, $\kappa(S_gF_n\cap F_n)\to \kappa(F_n)^2=1$ for each $n=1,\dots, N$ as $g\to\infty$.
Hence $\lim_{g\to\infty}\sum_{n=1}^\infty a_n|\lambda_n-1|(1-\kappa(S_gF_n\cap F_n))=0$,
as desired.
\qed

\comment

For that, it suffices to show that
$\langle U_{T_g^*}1,1\rangle\to 0$ as $g\to\infty$.
We let $f_k:=\frac{d\mu_k^*}{d\nu^*_k}$ and denote by $\langle .,\rangle_k$ the inner product in $L^2(Y^*,\mu_k^*)$.
Since $T_g^*=\big(\bigotimes_{1}^n S_g^*\big)\otimes\big(\bigotimes_{n+1}^\infty S_g^*\big)$ and $\sup_{n\in\Bbb N}\langle(\bigotimes_{k=n+1}^\infty U_{S_g^*})1,1\rangle_k\le 1$, it is enough to  prove that 
$\prod_{k=1}^n\big\langle U_{S_g^*}1,1\big\rangle_k\to 0$ as $g\to\infty$.
Since the dynamical system $(Y^*_k,\nu_k^*,S^*)$ is mixing, we obtain that
$$
\langle U_{S_g^*}1,1\rangle_k=\int_{Y^*}\sqrt{\frac{d\mu_k^*\circ U_{S_g^*}^{-1}}{d\mu^*_k}}\,
d\mu^*_k=\int_{Y^*}\sqrt{f_k\circ U_{S_g^*}^{-1} \cdot f_k}\,d\nu^*_k\to
\bigg(\int_{Y^*}\sqrt{f_k}d\mu^*_k\bigg)^2
$$
as $g\to\infty$.
It follows from Fact~1.18(iii) that for each $\omega\in Y^*$,
$$
f_k(\omega)=e^{-(\lambda_k-1)\nu_k(F_k)}\sum_{j=0}^\infty(\lambda_k)^j 1_{[F_k]_j}(\omega).
$$
Hence,
$$
\align
\int_{Y^*}\sqrt{f_k}d\nu_k^*&=e^{-0.5(\lambda_k-1)a_k}\sum_{j=0}^\infty\big(\sqrt{\lambda_k}\big)^j\nu_k^*([F_k]_j)\\
&=e^{-0.5(\lambda_k-1)a_k}\sum_{j=0}^\infty\big(\sqrt{\lambda_k}\big)^j e^{-a_k}\frac{a_k^j}{j!}\\
&=e^{-0.5(\lambda_k-1)a_k}e^{-a_k} e^{\sqrt{\lambda_k}a_k}\\
&=e^{-\frac{a_k}2(\sqrt{\lambda_k}-1)^2}.
\endalign
$$
Thus, we obtain that 
$
\langle U_{S_g^*}1,1\rangle_k\to e^{-a_k(\sqrt{\lambda_k}-1)^2}
$
as $g\to\infty$.
Therefore,
$$
\prod_{k=1}^n \langle U_{S_g^*}1,1\rangle_k\to e^{-\sum_{k=1}^n a_k(\sqrt{\lambda_k}-1)^2}.
\tag4-6
$$
It remains to note that $a_k(\sqrt{\lambda_k}-1)^2\sim\lambda_k^{-2}k^{-1}=\frac1{k4^{l_k}}$ and 
$\sum_{k=1}^n\frac1{k4^{l_k}}\to +\infty$ as $n\to \infty$, as desired.
\qed

\endcomment
\enddemo

\head 5. Type $III_\lambda$ ergodic Poisson suspensions for $\lambda\in(0,1)$
\endhead
Fix $\lambda\in(0,1)$.
In this section we
prove 
 the implications (1)$\Rightarrow$(3) of Theorems~A and B for  $K=III_\lambda$.
 Thus, we assume that $G$ does not have property (T).

Let $(Y,\goth Y,\kappa, S)$,  $\boldkey F=(F_n)_{n=1}^\infty$, $(g_k)_{k=1}^\infty$ and  $H$  denote the same objects as in~\S4.
Thus, $(\alpha_1)$--$(\alpha_5)$  and \thetag{4-1} from \S4 hold.
We now define $(a_n)_{n=1}^\infty$ and $(\lambda_n)_{n=1}^\infty$ in the following way:
\roster
\item"$(\alpha_6)'$" $\lambda_{2n-1}:=\lambda_{2n}^{-1}:=\lambda$ and 
$a_{2n-1}:=\lambda^{-1}a_{2n}:=\frac 1{n\log (n+1)}$ for each $n\in\Bbb N$.
\endroster

Denote by $(X,\goth B,\mu, T)$ the dynamical system associated with $(Y,\goth Y,\kappa,S)$,  $\boldkey F$, $(a_n)_{n=1}^\infty$ and $(\lambda_n)_{n=1}^\infty$ via the general construction from \S 2.
It follows from $(\alpha_3)$ and~$(\alpha_6)'$ that for each $g\in G$,
$$
\gather
\sum_{n\in\Bbb N}a_n|\lambda_n-1|\frac{\kappa
(S_gF_n\triangle F_n)}{\kappa(F_n)}<\infty\qquad\text{and}\tag5-1\\
\sum_{n=1}^\infty 
a_n\frac{
|1-\lambda_n^2|(1+\lambda_n^3)}{\lambda_n^2}
\frac{\kappa(S_gF_n\triangle F_n)}{\kappa(F_n)}<\infty.\tag5-2
\endgather
$$
Hence \thetag{5-1} and Lemma~2.1 yield that $T_g\in\text{Aut}_1(X,\mu)$  and $\chi(T_g)=0$ for each $g\in G$.
Therefore, the nonsingular Poisson suspension $(X^*,\goth B^*,\mu^*,T^*)$ is well defined by Fact~1.18(i).

\proclaim{Theorem~5.1} The dynamical system $(X^*,\goth B^*,\mu^*,T^*)$ is weakly mixing non-strongly ergodic   IDPFT of Krieger type $III_\lambda$.
Hence, $T^*$ is amenable in the Greenleaf sense. 
Moreover, $T^*$ is amenable if and only if $G$ is amenable.
\endproclaim
\demo{Proof}
Following the proof of Theorem~4.1, we obtain that $(X^*,\goth B^*,\mu^*,T^*)$ is IDPFT
and hence amenable in the Greenleaf sense. 
Moreover, we deduce from $(\alpha_1)$, $(\alpha_4)$ and  Propositions~2.3 and  2.4 that
$T^*$ is free and not strongly ergodic. 
 It follows from \thetag{5-2} and  Lemma~2.2    that
 $\bigg(\frac {d\mu}{d\mu\circ T_g^{-1}}\bigg)^2-1\in L^1(X,\mu)$ for each $g\in G$.
Since $c(\lambda^{-1})=\frac{c(\lambda)}\lambda$,  we deduce from $(\alpha_6)'$ that 
$$
c(\lambda_{2n})a_{2n}=c(\lambda_{2n-1})a_{2n-1}=\frac{c(\lambda)}{n\log(n+1)}\quad\text{ for every $n\in\Bbb N.$}
$$
Therefore for each $k\in\Bbb N$, in view of \thetag{4-1},
 $$
 \align
\int_X\Big(\Big(\frac{d\mu}{d\mu\circ T_{g_k}^{-1}}\Big)^2-1\Big)d\mu&=\sum_{n=1}^\infty a_nc(\lambda_n)\frac{\kappa
(S_gF_n\triangle F_n) }{2\kappa(F_n)}\\
&\le \frac{1}2 \sum_{n=1}^ka_nc(\lambda_n)+
\frac{c(\lambda)}2 \sum_{n=k+1}^\infty{\kappa
(S_{g_n}F_n\triangle F_n) }\\
&\le
 \frac14 \sum_{n=1}^k
\frac{c(\lambda)}{n\log(1+n/2)} +c(\lambda)\\
&\le c(\lambda)\log\log k +c(\lambda).
\endalign
$$
This implies that
$$
\sum_{k=1}^\infty \frac 1{(k\log k)^{1.5}}\,e^{\int_X\Big(\Big(\frac{d\mu}{d\mu\circ T_{g_k}^{-1}}\Big)^2-1\Big)d\mu}<+\infty.
$$
By Proposition~1.20 (see \thetag{1-4}),  
 $$
 \sum_{k=1}^\infty\bigg(\frac{d\mu^*\circ (T_{g_k}^*)^{-1}}{d\mu^*}(\omega)\bigg)^{\frac43}=+\infty\qquad\text{ at a.e. $\omega$.}
$$
On the other hand, $(Y^*,\nu_n^*, S^*)$ is mixing along $H$ because $(Y,\nu_n,S)$ is of 0-type along $H$ for each $n\in\Bbb N$ in view of Fact~1.21.
Therefore,  $(X^*,\mu^*,T^*)$ is weakly mixing in view of  Proposition~1.16. (Of course,
$
 \sum_{k=1}^\infty\frac{d\mu^*\circ (T_{g_k}^*)^{-1}}{d\mu^*}(\omega)=+\infty$ at a.e. $\omega$.)
Moreover,  the associated flow of $(X^*,\mu^*,T^*)$ coincides with the associated flow of $(X^*,\mu^*,\bigoplus_{n=1}^\infty  S^*)$ by Proposition~1.17.

{\sl Claim A.} The dynamical system $(X^*,\mu^*,\bigoplus_{n=1}^\infty  S^*)$ is of type $III_\lambda$.

We first show that $\lambda^{-1}$ is an essential value of  the Radon-Nikodym cocycle 
of $(X^*,\mu^*,\bigoplus_{n=1}^\infty  S^*)$.
Since  $\mu_{2k+1}(F_{2k+1})=\lambda_{2k+1}a_{2k+1}$, it follows from~$(\alpha_6)'$ that
for each  $n>0$,
$$
\sum_{k=n}^{+\infty}\mu_{2k+1}(F_{2k+1})=+\infty\quad\text{and}\quad\lim_{k\to\infty}\mu_{2k+1}(F_{2k+1})= 0.
$$
 Therefore,  there is $m_n>n$ such that 
$
\alpha_n:=\sum_{k=n}^{m_n-1}\mu_{2k+1}(F_{2k+1})\in (0.5,1).
$
Denote the probability space $((Y^*)^{2(m_n-n)}, \mu^*_{2n+1}\otimes\cdots\otimes\mu^*_{2m_n})$ by $(X^*_{2n+1,2m_n},\mu^*_{2n+1,2m_n})$.
We also let $\nu^*_{2n+1,2m_n}:=\nu^*_{2n+1}\otimes\cdots\otimes\nu^*_{2m_n}$.
Since 
$$
\mu_{2n-1}(F_{2n-1})-\nu_{2n-1}(F_{2n-1})+\mu_{2n}(F_{2n})-\nu_{2n}(F_{2n})=0,
$$
it follows from Fact~1.18(iii) that for a.e. $(\omega,\eta)\in Y^*\times Y^*$,
$$
\aligned
\frac{d(\mu_{2n-1}^*\otimes \mu_{2n}^*)}{d(\nu_{2n-1}^*\otimes\nu_{2n}^*)}(\omega,\eta)
&=\sum_{k=0}^\infty\lambda^k1_{[F_{2n-1}]_k}(\omega)
\sum_{k=0}^\infty\lambda^{-k}1_{[F_{2n}]_k}(\eta)\\
&=\sum_{k=-\infty}^\infty\lambda^k1_{B_{n,k}}(\omega,\eta),
\endaligned
\tag5-3
$$
where $1_{B_{n,k}}=\sum_{k=j-r} 1_{[F_{2n-1}]_j}1_{[F_{2n}]_r}$.
Hence the mapping
$$
\vartheta_n:X^*_{2n+1,2m_n}\ni \omega\mapsto\log_\lambda\frac{d\mu^*_{2n+1,2m_n}}{d\nu^*_{2n+1,2m_n}}(\omega)\in\Bbb Z
$$
is well defined.
Applying \thetag{5-3}, we obtain that
$$
\vartheta_n(\omega)=\sum_{k=n+1}^{m_n}(\omega(F_{2k-1})-\omega(F_{2k})),\quad\omega\in X^*_{2n+1,2m_n}.
$$
Thus, $\mu^*_{2n+1,2m_n}\circ\vartheta_n^{-1}$ is the distribution of the difference of two independent Poisson random variables $\omega\mapsto\sum_{k=n+1}^{m_n}\omega(F_{2k-1})$ and $\omega\mapsto\sum_{k=n+1}^{m_n}\omega(F_{2k})$, one with parameter $\alpha_n$, the other with
parameter $\lambda\alpha_n$.
In other words, $\mu^*_{2n+1,2m_n}\circ\vartheta_n^{-1}$ is the 
Skellam distribution with parameters $\alpha_n,\lambda\alpha_n$ (see \cite{AbSt}).
For $i=0,1$, let $\Delta_i:=\vartheta_n^{-1}(\{i\})\subset X^*_{2n+1,2m_n}$.
Then \footnote{Skellam (1946) and Prekopa (1953) represented the Skellam distribution using the modified Bessel function of the first kind. The result we are referring to is a direct consequence of this and standard facts on Bessel functions \cite{AbSt, pp. 374--378}.}
$$
\mu^*_{2n+1,2m_n}(\Delta_1)=
e^{-\alpha_n(1+\lambda)}\sum_{k=0}^\infty\frac{\alpha_n^{k+1}(\lambda\alpha_n)^k}{(k+1)!k!}>\frac{\alpha_n}{e^{\alpha_n(1+\lambda)}}>\frac1{16}
    $$
    and $\mu^*_{2n+1,2m_n}(\Delta_0)>\mu^*_{2n+1,2m_n}(\Delta_1)$.
Therefore,
$$
\frac{\nu^*_{2n+1,2m_n}(\Delta_0)}{\nu^*_{2n+1,2m_n}(\Delta_1)}=
\frac{\lambda\mu^*_{2n+1,2m_n}(\Delta_0)}{\mu^*_{2n+1,2m_n}(\Delta_1)}\ge \lambda.
$$
Since the action $\bigoplus_{j=2n+1}^{2m_n} S^*$ of $G^{2(m_n-n)}$ on $\big((Y^*)^{2(m_n-n)},\nu^*_{2n+1,2m_n}\big)$ is measure preserving and ergodic,
there exist  a subset $\Delta'_1\subset\Delta_1$  and a transformation
$Q$ from the full group $\big[\bigoplus_{j=2n+1}^{2m_n} S^*,\nu^*_{2n+1,2m_n}\big]$ such that $Q\Delta'_1\subset\Delta_{0}$ and $\nu^*_{2n+1,2m_n}(\Delta'_1)=\lambda\nu^*_{2n+1,2m_n}(\Delta_1)$.
Then
$$
\mu^*_{2n+1,2m_n}(\Delta_1')=\lambda\nu^*_{2n+1,2m_n}(\Delta_1')=\lambda^2\nu^*_{2n+1,2m_n}(\Delta_1)=\lambda\mu^*_{2n+1,2m_n}(\Delta_1)>\frac\lambda{16}.
$$
If $\omega\in \Delta_1'$ then  $Q\omega\in \Delta_0$ and hence
$\frac{d\mu^*_{2n+1,2m_n}}{d\nu^*_{2n+1,2m_n}}(Q\omega)=1$ and
$\frac{d\mu^*_{2n+1,2m_n}}{d\nu^*_{2n+1,2m_n}}(\omega)=\lambda$ by the definition of 
$\Delta_0$ and $\Delta_1$.
Therefore
$$ 
\frac{d\mu^*_{2n+1,2m_n}\circ Q}{d\mu^*_{2n+1,2m_n}}(\omega)=
\frac{d\mu^*_{2n+1,2m_n}}{d\nu^*_{2n+1,2m_n}}(Q\omega)\frac{d\nu^*_{2n+1,2m_n}}{d\mu^*_{2n+1,2m_n}}(\omega)=\lambda^{-1}.
$$
Now take a Borel subset  $C$ of $(Y^*)^{2n}$.
Then
\roster
\item"---"
$[C\times\Delta'_1]_{2m_n}\subset [C]_{2n}$,
\item"---"
$\mu^*([C\times\Delta'_1]_{2m_n})=\mu^*([C]_{2n})\mu^*_{2n+1,2m_n}(\Delta_1')>\frac\lambda {16}\mu^*([C]_{2n})$,
\item"---"
$I\times Q\times I\in\big[\bigoplus_{n=1}^\infty  S^*,\mu^*\big]$ and 
$
\frac{d\mu^*\circ (I\times Q\times I)}{d\mu^*}(\omega)=\lambda^{-1}
$
 for $\mu^*$-a.e. $\omega\in[C\times\Delta'_1]_{2m_n}$.
\endroster
Since the family $\{[C]_{2n}\mid C\subset (Y^*)^{2n}, n\in\Bbb N\}$ is dense
in the entire Borel $\sigma$-algebra on $X^*$, it follows  that $\lambda^{-1}\in r\big(\bigoplus_{n=1}^\infty  S^*,\mu^*\big)$ by Fact~1.7.
On the other hand,  the Radon-Nikodym cocycle of the system $(X^*,\mu^*,\bigoplus_{n=1}^\infty  S^*)$ 
 takes its values in the subgroup $\{\lambda^n\mid n\in\Bbb Z\}$ of $\Bbb R_+^*$.
  It follows that $r\big(\bigoplus_{n=1}^\infty  S^*,\mu^*\big)\subset\{\lambda^n\mid n\in\Bbb Z\}$.
 Hence $r\big(\bigoplus_{n=1}^\infty  S^*,\mu^*\big)=\{\lambda^n\mid n\in\Bbb Z\}$.
 Therefore,
 $(X^*,\mu^*,\bigoplus_{n=1}^\infty  S^*)$ is of type $III_\lambda$, as desired.
 Thus, Claim~A is proved.

 It follows from Claim A and Proposition~1.17 that  $(X^*,\mu^*,T^*)$ is of type $III_\lambda$.
\qed
\enddemo

\remark{Remark \rom{5.2}}  Arguing as in Remark~4.2, one can choose the parameters $(a_n)_{n=1}^\infty$ and $(\lambda_n)_{n=1}^\infty$ in such a way that  the Poisson suspension $T$ of the associated action $T$ (we mean the general construction of \S2) is of infinite ergodic index indeed.
We leave details to the reader.
\endremark


\proclaim{Theorem 5.3}
If $G$ is Haagerup then there is a dynamical system $(X,\goth B,\mu, T)$ such that the Poisson suspension $( X^*,\goth B^*,\mu^*,  T^*)$ of $(X,\goth B,\mu, T)$  is 
of 0-type (in addition to the other properties of $( X^*,\goth B^*,\mu^*,  T^*)$ listed in Theorem~5.1).
\endproclaim

\demo{Proof}
Since $G$ is Haagerup,
 there is a free 0-type $G$-action $S=(S_g)_{g\in G}$ on an infinite $\sigma$-finite
standard measure space $(Y,\goth Y,\kappa)$ such that $(\alpha_1)$--$(\alpha_5)$, \thetag{4-1} and $(\alpha_6)'$ hold.

We fix an increasing sequence $H_1\subset H_2\subset \cdots$ of compact subsets in $G$
with $\bigcup_{n=1}^\infty H_n=G$.
Then there are a  sequence  $K_1\subset K_2\subset \cdots$ of compact subsets in $G$
with $\bigcup_{n=1}^\infty K_n=G$
and sequence of positive reals $n_1<n_2<\cdots$ such that for each $l\in\Bbb N$,
\roster
\item"$(\alpha_7)$"
$\sup_{g\in K_l}\kappa(S_gF_{n_{l}}\triangle F_{n_l})<\frac1{l^2}$ and
\item"$(\alpha_8)$"
$\sup_{1\le j\le l}\sup_{g\not\in K_{l+1}}(1-\kappa(S_gF_{n_j}\cap  F_{n_j}))<\frac1{l}$.
\endroster
These sequences can be constructed inductively.
If we have already determined $(n_j)_{j=1}^{l-1}$  and $(K_j)_{j=1}^l$ then there exists $n_l$ such that $(\alpha_7)$ holds.
This follows from the fact that $(F_n)_{n=1}^\infty$ is $S$-F{\o}lner and $\kappa(F_n)=1$ for each $n\in\Bbb N$.
Then we select  a compact subset $K_{l+1}\subset G$ large so that
$K_{l+1}\supset K_l\cup H_l$  and $(\alpha_8)$ holds.
This follows from the fact that $S$ is of 0-type and  $\kappa(F_n)=1$ for each $n\in\Bbb N$.
Repeating these 2 steps infinitely many times we determine the entire sequences
$(n_j)_{j=1}^{\infty}$  and $(K_j)_{j=1}^{\infty}$.

To simplify notations, we rename the subsequence $(F_{n_l})_{l=1}^\infty$ into $(F_n)_{n=1}^\infty$.
Thus, from now on we may assume without loss of generality that $(\alpha_1)$--$(\alpha_5)$, \thetag{4-1} and~$(\alpha_6)'$ hold and, in addition, for each $n\in\Bbb N$,
\roster
\item"$(\alpha_7)'$"
$\sup_{g\in K_n}\kappa(S_gF_{n}\triangle F_n)<\frac1{n^2}$ and
\item"$(\alpha_8)'$"
$\sup_{1\le j\le n}\sup_{g\not\in K_{n+1}}(1-\kappa(S_gF_{j}\cap  F_{j}))<\frac1n$.
\endroster

Repeating the proof of Theorem~5.1 almost literally, we construct a system $(X,\goth B,\mu,T)$  such that the Poisson suspension  $(X^*,\goth B^*,\mu^*,T^*)$ of 
 $(X,\goth B,\mu,T)$ is well defined and possesses  all the properties listed in 
 Theorem~5.1.
 It remains to prove that  $( X^*,\goth B^*,\mu^*,  T^*)$ is of 0-type.
For that, we first show that
 $$
\lim_{g\to\infty} \sum_{n=1}^\infty \frac 1{n\log(n+1)}\kappa(S_gF_n\triangle F_n)= 0.\tag5-4
 $$ 
 Given $\epsilon>0$, we  select $N>0$ such that $\sum_{n>N}\frac 1 {n^2}<\epsilon$.
 Let $g\not\in K_N$.
 Then there is $l\ge N$ such that $g\in K_{l+1}\setminus K_l$.
 It follows from $(\alpha_7)'$ that
$$
\sum_{n>l} \frac 1{n\log(n+1)}\kappa(S_gF_n\triangle F_n)<\sum_{n>l} \frac 1{n^2}<\epsilon.
$$
On the other hand, it follows from $(\alpha_8)'$ that
 $$
\sum_{n=1}^l \frac 1{n\log(n+1)}\kappa(S_gF_n\triangle F_n)<
\frac1{l-1}\sum_{n=1}^{l-1}\frac 1{n\log(n+1)} + \frac{2}{l\log(l+1)}
$$
 Since the righthand side of this inequality goes is less than $\epsilon$ whenever $N$ is large enough, \thetag{5-4} follows.
We now deduce from  \cite{DaKoRo1, Proposition~6.9}, $(a_6)'$ and  \thetag{5-4} 
that $T^*$ is of zero type.
\qed
\enddemo

\head 6. Type $III_1$ ergodic Poisson suspensions 
\endhead
In this section we
prove 
 the implications (1)$\Rightarrow$(3) of Theorems~A and B for  $K=III_1$.
 Thus, we assume that $G$ does not have property (T).
Fix $\lambda_1,\lambda_2\in(0,1)$ such that $\log\lambda_1$ and $\log\lambda_2$ are rationally independent.
 The following theorem  follows from Theorem~5.1 (as above, we use here the notation from \S2).
 For $n\in\Bbb N$, let 
 $$
 \align
 \lambda_{4n-3}&:=\lambda_{4n-2}^{-1}:=\lambda_1, \\
\lambda_{4n-1}&:=\lambda_{4n}^{-1}:=\lambda_2,\\
a_{4n-3}&:=\lambda^{-1}_1a_{4n-2}:=\frac 1{n\log (n+1)} \quad\text{
and }\\
a_{4n-1}&:=\lambda_2^{-1}a_{4n}:=\frac 1{n\log (n+1)}.
\endalign
$$
Denote by $(X,\goth B,\mu, T)$ the dynamical system associated with $(Y,\goth Y,\kappa,S)$,  $\boldkey F$, $(a_n)_{n=1}^\infty$ and $(\lambda_n)_{n=1}^\infty$ via the general construction from \S 2.

\proclaim{Theorem 6.1} The Poisson suspension  $(X^*,\goth B^*,\mu^*,T^*)$ 
of $(X,\goth B,\mu, T)$ is a well defined nonsingular dynamical system.
It is free, non-strongly ergodic,  IDPFT,  of Krieger type $III_1$ and of infinite ergodic index.
Hence  $T^*$ is weakly mixing and amenable in the Greenleaf sense.
Moreover, $T^*$ is amenable if and only if $G$ is amenable.
\endproclaim

\demo{Idea of the proof} 
Follow the proof of  Theorem~5.1 to show that
the associated flow of  $(X^*,\goth B^*,\mu^*,T^*)$ is isomorphic to the associated flow
of the the direct sum  $(X^*,\goth B^*,\mu^*,\bigoplus_{n=1}^\infty S^*)$ by Proposition~1.17.
However, the latter system splits into the direct sum (do not confuse with the direct product)
of the following two systems:$\big(\bigotimes_{j\in I} (Y^*,\mu_j^*),\bigoplus_{j\in I} S^*\big)$ and  $\big(\bigotimes_{j\in J} (Y^*,\mu_j^*),\bigoplus_{j\in J} S^*\big)$,
where 
$$
\align
I&:=\{4n-2\mid n\in\Bbb N\}\cup\{4n-3\mid n\in\Bbb N\}\text{ and }\\
J&:=\{4n-1\mid n\in\Bbb N\}\cup\{4n\mid n\in\Bbb N\}.
\endalign
$$ 
These two dynamical systems  are of type $III_{\lambda_1}$ and $III_{\lambda_2}$ respectively
according to Claim~A from the proof of Theorem~5.1.
Since $\log\lambda_1$ and $\log\lambda_2$ are rationally independent,
the associated flow of $(X^*,\goth B^*,\mu^*,\bigoplus_{n=1}^\infty S^*)$ is trivial (acting on a singleton).
Hence $(X^*,\goth B^*,\mu^*,T^*)$
  is of type $III_1$.
  The other properties of this system can be established in the very same way as in the proof
  of Theorem~5.1.
\qed
\enddemo

In a similar way one can prove the following assertion.

\proclaim{Theorem 6.2}
If $G$ is Haagerup then there is a dynamical system $(X,\goth B,\mu, T)$ such that the Poisson suspension $( X^*,\goth B^*,\mu^*,  T^*)$ of $(X,\goth B,\mu, T)$ is a well defined nonsingular dynamical system.
This system is of 0-type and possesses all the properties of $( X^*,\goth B^*,\mu^*,  T^*)$ listed in Theorem~6.1.
\endproclaim

\head 7. Applications to nonsingular Bernoulli actions of locally compact amenable groups 
\endhead

Let $G$ be a  non-compact locally compact second countable group.
Given an infinite $\sigma$-finite standard measure space $(X,\goth B,\mu)$, we let
(as in \cite{DaKoRo1})
$$
\text{\rom{Aut}}_2(X,\mu):=\bigg\{S\in\text{\rom{Aut}}(X,\mu)\,\bigg|\,\sqrt{\frac{d\mu\circ S}{d\mu}}-1\in L^2(X,\mu)\bigg\}.
$$
As was shown in \cite{DaKoRo1}, given a transformation $S\in\text{\rom{Aut}}(X,\mu)$,
 the nonsingular Poisson suspension $S^*$ of $S$ is well defined (as an element of $\text{\rom{Aut}}(X^*,\mu^*)$) if and only if $S\in\text{\rom{Aut}}_2(X,\mu)$.
 In particular,
  Aut$_1(X,\mu)\subset\text{Aut}_2(X,\mu)$.

\definition{Definition 7.1} Let $T=(T_g)_{g\in G}$ be a free totally dissipative nonsingular 
$G$-action on a $\sigma$-finite standard measure space $(X,\goth B,\mu)$
such that $T_g\in \text{\rom{Aut}}_2(X,\mu)$ for each $g\in G$.
Then we call the Poisson suspension $T^*:=(T^*_g)_{g\in G}$ of $T$  {\it the nonsingular Bernoulli $G$-action} over the base $(X,\goth B,\mu,T)$.
\enddefinition

If $G$ is  discrete countable then each nonsingular Bernoulli $G$-action according to Definition~7.1 is nonsingular Bernoulli in the usual sense.
However, we do not know whether each nonsingular Bernoulli $G$-action in the usual sense
is isomorphic to a Bernoulli $G$-action in the sense of Definition~7.1.

For an arbitrary $G$, if $T$ preserves $\mu$ then $T^*$ is a probability preserving Bernoulli $G$-action  in the sense of \cite{OrWe}.

Let $\lambda_G$ denote  a left Haar measure on $G$.
Since $T$ is totally dissipative and free, there is a  Borel subset $B\subset X$ which meets a.e. $T$-orbit exactly once.
Hence there exist a $\sigma$-finite Borel measure $\kappa$ on $B$ and  a Borel isomorphism (mod 0) of $X$ onto $B\times G$ such that $T$ corresponds to $G$-action by the left translations along the second coordinate of $B\times G$.
We identify  $X$ and $B\times G$ via this isomorphism.
Choose a $\sigma$-finite Borel measure $\kappa$ on $B$ which has the same collection of subsets of 0 measure as the  projection of $\mu$ to $B$ has.
Then the disintegration $\mu=\int_B\delta_b\otimes\mu_b\,d\kappa(b)$ of $\mu$ with respect to $\kappa$ is well defined.
Of course, $\mu_b\sim\lambda_G$ at a.e. $b\in B$.
Therefore, there is a measurable function $F:B\times G\to\Bbb R^*_+$ such that
$F(b,g)=\frac{d\mu_b}{d\lambda_G}(g)$ almost everywhere.
The condition $T_h\in \text{\rom{Aut}}_2(X,\mu)$ can be rewritten now as
$$
\int_{B}\int_G(\sqrt{F(b,g)}-\sqrt{F(b,hg)})^2\,d\lambda_G(g)d\kappa(b)<\infty
$$
for each $h\in G$.
In particular, letting $c_{b,h}(g):=\sqrt{F(b,g)}-\sqrt{F(b,hg)}$, we obtain that $c_{b,h}\in L^2(G,\lambda_G)$ for each $b\in B$ and $h\in H$.
Moreover, the function  $B\ni b\mapsto\|c_{b,h}\|_2$ belongs to $L^2(B,\kappa)$.

\definition{Definition 7.2}
If there is a measurable function $d:B\ni b\mapsto d_b\in L^2(G,\lambda_G)$ such that
$c_{b,h}(g)=d_b(g)-d_b(hg)$ at a.e. $g\in G$ and each $h\in H$ then we say that  $(X,\goth B,\mu, T)$ is {\it tame}.
\enddefinition

If  $(X,\goth B,\mu, T)$ is tame then we can find a countable partition $B=\bigsqcup_{j}B_j$ of $B$ into measurable subsets such that the function 
$B_j\ni b\mapsto\|d_{b}\|_2$ belongs to $L^2(B_j,\kappa\restriction B_j)$ for each $j$.
Then we get a countable partition $\bigsqcup_{j}(B_j\times G)$ of $X$ into $T$-invariant
subsets  furnished with $T$-invariant measures $\nu_j$ such that $\sqrt{\frac{d\mu}{d\nu_j}}-1\in L^2(B_j\times G,\nu_j)$ for each $j$.
It follows that $((B_j\times G)^*,(\mu\restriction(B_j\times G))^*, T^*)$ is  nonsingular Bernoulli of type $II_1$.
Thus, if $(X,\goth B,\mu, T)$ is tame then the corresponding Bernoulli $G$-action over $(X,\goth B,\mu, T)$ is IDPFT.

\proclaim{Corollary 7.3} Let $G$ be amenable non-compact locally compact second countable group.
For each $K\in\{III_\lambda\mid 0\le\lambda\le 1\}\cup\{II_\infty\}$,  there is a   tame Bernoulli  free nonsingular $G$-action  $T$ of infinite ergodic index and of Krieger type $K$.
Hence $T$ is weakly mixing, IDPFT and amenable in the Greenleaf sense.
In particular, $T$ is amenable if and only if $G$ is amenable.
\endproclaim

\demo{Proof}
Since $G$ is amenable then there is a left F{\o}lner sequence $(F_n)_{n=1}^\infty$ on $G$.
Let $(Y,\kappa):=(G,\lambda_G)$.
Let $S$ denote  the  $G$-action on $Y$ by left translations and
let
$\boldkey F:=(F_n)_{n=1}^\infty$.
Then the Koopman representation $U_S$ is the left regular representation of $G$.
Hence $U_S$ is mixing, i.e.  $S$ is of 0-type.
Moreover, $\boldkey F$ is $S$-F{\o}lner.
It remains to apply the constructions from Theorems~3.5, 4.3, 5.3 and 6.2 (depending on the value of $K$) to the quadruple  $(Y,\kappa,S,\boldkey F)$.
\qed

\comment

Utilizing the constructions from the proof of Theorem~3.1 for the system $(Y,\kappa,S)$ and $\boldkey F$, we obtain a totally dissipative dynamical system
$(\widetilde X,\widetilde{\goth B},\widetilde\mu,\widetilde T)$ which is tame.
Hence the Poisson suspension of this system is a weakly mixing Bernoulli IDPFT $G$-action of type $II_\infty$ according to Theorem~3.1 and the remark above~Corollary~7.3.

As for the type $III_0$, we follow the proof of Theorem~4.1.
However, in the (amenable) case under consideration here, the argument is simpler as we start from $(Y,\kappa,S)$ as in the $II_\infty$-case as above and do not need the Moore-Hill construction over an auxiliary Gaussian $G$-action.
Let $R$ denote the right $G$-action on $(Y,\kappa)$.
Since $G$ is unimodular, $R$ preserves $\kappa$.
Moreover, $R$ is ergodic and commutes with $S$.
Thus, $R$ plays now  the role of $\Gamma$ in the proof of Theorem~4.1.
Repeating the rest of the proof of Theorem~4.1 almost verbally, we deduce that the Poisson suspension of a tame totally dissipative dynamical system $(X,\goth B,\mu,T)$ constructed
with usage of  $(Y,\kappa,S)$ is a  weakly mixing Bernoulli $G$-action of type $III_0$.

The remaining types $III_\lambda$, $0<\lambda\le 1$ are treated almost verbally as in Theorems~5.1 and 6.1 but with $(Y,\kappa,S)$ and $\boldkey F$ as above.

\endcomment

\enddemo

\head 8. Interplay between nonsingular Poisson and nonsingular Gaussian actions
\endhead
Given a real separable Hilbert space $\Cal H$, 
we denote by $\Cal A$ the group of affine transformations of $\Cal H$, i.e.
$\Cal A=\Cal H\rtimes\Cal O$, where $\Cal O$ is the group of orthogonal transformations of $\Cal H$.
A transformation $(h,O)\in \Cal A$ acts on a vector $v\in\Cal H$ by the formula
$(h,O)v:=h+Ov$.
It was shown in \cite{ArIsMa}  (see also \cite{DaLe, Section~3}) that given a continuous homomorphism $\alpha:G\to\Cal A$, one can construct a nonsingular $G$-action, called {\it the nonsingular Gaussian} action generated by $\alpha$.

Suppose now that we have  
an infinite $\sigma$-finite standard measure space $(Y,\goth Y,\nu)$. 
Let $S=(S_g)_{g\in G}$ be  a nonsingular $G$-action on $Y$ such that $S_g\in\text{Aut}_2(Y,\nu) $ for all $g\in G$.
Then the Poisson suspension $S^*=(S_g^*)_{g\in G}$ is well defined as a nonsingular $G$-action on $(X^*,\goth Y^*,\nu^*)$ \cite{DaKoRo1}.
Denote by $\Cal A$ the group of affine transformations of $L^2(Y,\nu)$.
Then a continuous homomorphism 
$$
A_S:G\ni g\mapsto A_S(g)\in \Cal A
$$
 is well defined 
by the formula:
$
A_S(g):=\Big(2\Big(\sqrt{\frac{d\nu\circ S_g^{-1}}{d\nu}} -1 \Big),U_S(g)\Big)\in \Cal A$.
We call $A_S$ {\it the affine Koopman representation} of $G$ associated with $S$ \cite{DaKoRo1}.
Denote by 
 $H=(H_g)_{g\in G}$ the nonsingular Gaussian $G$-action generated by $A_S$.
It was explained in \cite{DaLe, Remark~3.6} that  $H$ and $S^*$ are spectrally identical, i.e. the unitary Koopman representations $U_H$ and $U_{S^*}$ of $G$ are unitarily equivalent.
Hence it is of interest to compare non-spectral dynamical properties of $H$ and $S^*$.

\proclaim{Proposition 8.1} Let $G$  have the Haagerup property.
Let $S$ be a 0-type nonsingular $G$-action on the standard $\sigma$-finite measure space $(Y,\goth Y,\nu)$ such that $S_g\in\text{\rom{Aut}}_2(Y,\nu)$ for each $g\in G$.
 Suppose that  $Y=\bigsqcup_{n=1}^\infty Y_n$ for some $S$-invariant subsets $Y_n\in\goth Y$ of infinite measure $\nu$
 and
 there is an $S$-invariant $(\nu\restriction Y_n)$-equivalent 
$\sigma$-finite measure $\xi_n$ on $Y_n$ with $\sqrt{\frac{d\xi_n}{d\nu}}-1\in L^2(Y_n,\nu)$ for each $n\in\Bbb N$.
If  the Poisson suspension $(Y^*,\nu^*,S^*)$ does not admit an invariant equivalent probability measure then 
 $H$ is either totally dissipative or weakly mixing of type $III_1$.
 In particular, if $G$ is non-amenable then $H$ is
 weakly mixing of type $III_1$.
\endproclaim
\demo{Proof}  
We first note that $L^2(Y,\nu)=\bigoplus_{n=1}^\infty L^2(Y_n,\nu)$ and that $L^2(Y_n,\nu)$ is an invariant subspace for $U_S$ for 
each $n\in\Bbb N$.
Let 
$$
c_g:=2\Bigg(\sqrt{\frac{d\nu\circ S_g^{-1}}{d\nu}} -1\Bigg)\in L^2(Y,\nu)\qquad\text{for each $g\in G$.}
$$
Consider $ L^2(Y,\nu)$ as a $G$-module, where $G$ acts via $U_S$.
Then the mapping 
$$
c:G\ni g\mapsto c_g\in L^2(Y,\nu)
$$ is a 1-cocycle of $G$ with coefficients in $L^2(Y,\nu)$.
It follows from the condition of the proposition and \cite{DaKoRo1, Proposition~6.4} that for each $n>0$, 
the cocycle $G\ni g\mapsto 1_{Y_n}c_g\in L^2(Y_n,\nu)$ is a coboundary.
Hence there exists a function $a_n\in L^2(Y_n,\nu)$ such that
$1_{Y_n}c_g= a_n-U_S(g)a_n$ for each $g\in G$.
We now set  
  $$
  \Cal H_n^0:=\bigoplus_{m>n}L^2(Y_m,\nu)\qquad\text{and}\qquad 
\Cal H_n:=a_1\oplus\cdots \oplus a_n\oplus\Cal H_n^0.
$$
Then $\Cal H_1^0\supset \Cal H_2^0\supset\cdots$, $\bigcap_{n=1}^\infty\Cal H^0=\{0\}$ and
 $$
 A_S(g)\Cal H_n=\bigoplus_{m=1}^n(U_S(g)a_m+1_{Y_m}c_g)\oplus U_S(g)\Cal H^0_n=\bigoplus_{m=1}^n a_n\oplus\Cal H_n^0=\Cal H_n
 $$
for each $g\in G$ and $n\in\Bbb N$.
Thus, $\Cal H_n$ is invariant under $A_S$.
Hence, by \cite{ArIsMa, Proposition~2.10}, $A_S$  is {\it evanescent} according to \cite{ArIsMa, Definition~2.6}.
Since $S$ is of 0-type, $U_S$ is mixing.
Since $(Y^*,\nu^*,S^*)$ does not admit an invariant equivalent probability measure then 
$c$ is not a coboundary \cite{DaKoRo1, Proposition~6.4}.
Therefore  \cite{ArIsMa, Theorem~D} yields that $H$ is either totally dissipative or weakly mixing and 
 of Krieger type $III_1$.
 The first claim of the proposition is proved.
 
 Suppose now that $G$ is non-amenable.
 We have to show  that $H$ is not totally dissipative.
 It follows from the condition of the proposition that  $(Y^*,\nu^*,S^*)$ is IDPFT.
 Hence by Proposition~1.25, $S^*$ is amenable in the Greenleaf sense.
 Since the amenability in the Greenleaf sense is an invariant for the unitary equivalence,
 it follows that  $H$ is also amenable in the Greenleaf sense. 
We deduce from this fact (as $G$ is nonamenable) that  $H$ is not amenable in view of Fact~1.24 (ii), (iv) and (v).
Since each totally dissipative $G$-action is amenable, $H$ is not totally dissipative.
 \qed
\enddemo

\proclaim{Corollary 8.2} Let  $G$ be  non-amenable and have the Haagerup property.
Then for all nonsingular Poisson $G$-actions  (of any Krieger type) constructed in Theorem~B, the corresponding 
nonsingular Gaussian $G$-actions are all weakly mixing and of Krieger type $III_1$.
\endproclaim

\comment

\head 9. Proof of Corollary E
\endhead

Let $G$ be non-amenable and non-(T).
Take $K\in\{III_\lambda\mid \lambda\in[0,1]\}\sqcup\{II_\infty\}$.
Let $T^*$ be the corresponding Poisson $G$-action
 whose existence is claimed in Theorem~A.
 We note that $T^*$ appears as the Poisson suspension of a system $(X,\mu,T)=\bigsqcup_{n\in\Bbb N}(Y,\mu_n,S)$ such that 
 \roster
 \item"---"
 there  is an infinite   $\sigma$-finite  $S$-invariant measure $\nu_n\sim\mu_n$ on $Y$ with $\frac{d\mu_n}{d\nu_n}-1\in L^1(Y,\nu_n)$ for each $n\in\Bbb N$ and 
 \item"---"
 there is an $S$-F{\o}lner sequence in $Y$.
 \endroster
It was explained in the beginning of Section~4 that we can also assume without loss of generality that 
\roster 
\item"---" there exists a countable ergodic group $\Gamma$ of $\nu_n$-preserving transformations (for each $n$) commuting with $S$.\footnote{Though this fact was used only in Section~4.}
\endroster
By Theorem~A, $T^*$ is weakly mixing, non-strongly ergodic  and non-amenable (because $G$ is non-amenable).

It remains to show that $T^*$ is free.
Since $(X^*,\mu^*,T^*)=\bigotimes_{n=1}^\infty (Y^*,\mu_n^*,S^*)$
and  $\mu_1^*\sim\nu_1^*$ by Fact~1.18(ii),
it suffices to show that the action
$S^*$ on $(Y^*,\nu_1^*)$ is free.
Let $\Gamma^*=\{\gamma^*\mid\gamma\in \Gamma\}$.
Then $\Gamma^*$ preserves $\nu_1^*$.
Since $(Y,\nu_1,\Gamma)$ is ergodic, the dynamical system $(Y^*,\nu_1^*,\Gamma^*)$ is ergodic too \cite{Ro}.
 Moreover, $\Gamma^*$ commutes with $S^*$.
 Denote by $\Cal G$ the space of all closed subgroups of $G$ and endow $\Cal G$ with the Fell topology \cite{Fe}.
 Then $\Cal G$ is a compact metric space \cite{Fe}.
 Given $\omega\in Y^*$, let $G_\omega:=\{g\in G\mid S^*_g\omega=\omega\}$
 stand for  the stability group of $S^*$ at $\omega$.
 Then $G_\omega\in\Cal G$ at each $\omega$ \cite{AuMo, I, Proposition~3.7} and the mapping
 $
 \eta:Y^*\ni\omega\mapsto G_\omega\in \Cal G
 $
is  Borel \cite{AuMo, II, Proposition~2.3}.
Since $\Gamma^*$ commutes with $S^*$, it is straightforward to verify that  $\eta$ is invariant under $\Gamma^*$.
Hence, $\eta$ is constant, i.e. there is a subgroup $H\in\Cal G$ such that $G_\omega=H$ at a.e. $\omega\in Y^*$.
Therefore $S^*_g=\text{Id}$ for each $g\in H$.
Since $S^*$ is effective, $H$ is trivial, i.e. $S^*$  (and hence $T^*$) is free, as desired.

\endcomment

\enddocument